\documentclass[12pt,a4paper]{article}
\pdfoutput=1
\usepackage{amsmath,amssymb,amsthm}
\usepackage{apropos,graphicx,pict2e}
\usepackage[left=20mm,right=20mm,top=25mm,bottom=25mm]{geometry}

\usepackage[numbers,sort&compress]{natbib}

\usepackage{array}
\usepackage{algorithm}
\usepackage{algpseudocode}

\title{Automatically Stable Discontinuous Petrov-Galerkin Methods for Stationary Transport Problems: Quasi-Optimal Test Space Norm}

\author{Antti H. Niemi\thanks{antti.niemi@kaust.edu.sa}, 
Nathaniel O. Collier\thanks{nathaniel.collier@kaust.edu.sa}, 
Victor M. Calo\thanks{victor.calo@kaust.edu.sa} \\ \\
King Abdullah University of Science and Technology (KAUST) \\
Thuwal 23955-6900, Kingdom of Saudi Arabia
}
\date{}
\begin{document}

\maketitle

\begin{abstract}
We investigate the application of the discontinuous Petrov-Galerkin (DPG) finite element framework to stationary convection-diffusion problems. In particular, we demonstrate how the quasi-optimal test space norm can be utilized to improve the robustness of the DPG method with respect to vanishing diffusion. We numerically compare coarse-mesh accuracy of the approximation when using the quasi-optimal norm, the standard norm, and the weighted norm. Our results show that the quasi-optimal norm leads to more accurate results on three benchmark problems in two spatial dimensions. We address the problems associated to the resolution of the optimal test functions with respect to the quasi-optimal norm by studying their convergence numerically. In order to facilitate understanding of the method, we also include a detailed explanation of the methodology from the algorithmic point of view.
\end{abstract}


\section{Introduction}
The purpose of this paper is to study the application of the discontinuous Petrov-Galerkin finite element methodology to convection-diffusion problems of the form
\begin{equation} \label{eqn:the equation}
\begin{aligned}
- \epsilon \Delta u + \va \cdot \vnabla u &= f \quad \text{in $\Omega$} \\
u &= g \quad \text{on $\partial \Omega$}
\end{aligned}
\end{equation}
where $\Omega$ is a bounded two-dimensional domain, $\va$ is the flow velocity, $\epsilon$ is the diffusion parameter and $f$ represents a source term. The function $g$ determines the values of the solution $u$ on the boundary of the domain denoted by $\partial\Omega$. This equation is a basic model for various transport processes in nature and engineering. Typical applications include the movement of substances (dissolved nutrients, toxins) in the environment. Problem \eqref{eqn:the equation} may also be viewed as a stepping stone for designing numerical methods for more complex problems such as the Navier-Stokes equations for compressible and incompressible flows.

Construction of numerical solution schemes for \eqref{eqn:the equation} that can approximate convection dominated flows is known to be an inherently difficult task. Many scientists and engineers have addressed the problem and designed algorithms based on finite difference, finite volume, finite element and spectral approaches. In the context of the finite element method, the standard Ritz-Bubnov-Galerkin formulation is known to be deficient and various alternative formulations have been proposed in the literature, see for instance \cite{HeiHuyZieMitIJNME1977,BroHugCMAME1982,JohCMAME1992,DemOdeJCOMPUTPHYS1986,DemOdeCMAME1986,Jiang1993,Guermond1999,BreFraHugRusCMAME1997,CodCMAME1998,OnaCMAME1998,Brezzi2000,Guermond2004,BufHugSanSJNA2006,HugEvaSanCMAME2009}.

The success of the traditional Galerkin finite element method in most structural problems is based on its best approximation property. This means that the difference between the finite element solution and the exact solution is minimized with respect to a norm induced by the strain energy product. Numerical issues arise when the best approximation property is lost.  This happens, for instance, when the standard Galerkin method is applied to convective transport problems. In these problems, the system matrix associated to convection is not symmetric and numerical solutions tend to show spurious, non-physical oscillations unless the finite element mesh is heavily refined. It should be noted that the accuracy of the finite element method may also degrade in problems with symmetric system matrices. Thin-body problems of solid mechanics and wave propagation problems constitute typical examples. In these problems, the standard Galerkin method leads to non-physical solutions suffering from the phenomena of numerical locking (small thickness) and numerical dispersion (high wave number) unless an over-refined mesh is used.

To save computer resources, alternative formulations which produce reasonable approximations regardless of the mesh size have been developed. In the context of convective transport problems, a well-known method is the Streamline Upwind Petrov-Galerkin (SUPG) formulation, see \cite{BroHugCMAME1982,JohNavPitCMAME1984}. The method belongs to a larger class of stabilized finite element methods reviewed, for example, in \cite{CodCMAME1998,FraHauMasCMAME2006}. The stabilized finite element methods are geared mainly towards the classical $h$-version of the finite element method although some higher order formulations have been studied in \cite{CodCIMNE1993,FraFreHugCMAME1992,HugCotBazCMAME2005,BazCalCotEvaHugLipScoSedCMAME2010}.

A new, general finite element framework has been developed by Demkowicz and Gopalakrishnan in \cite{DemGopCMAME2010,DemGopNUMPDE2010,DemGopIACM2011}. This variational framework is of discontinuous Petrov-Galerkin (DPG) type and can be utilized, for example, to symmetrize non-symmetric problems by using the concept of optimal test functions similar to the early methodology developed by Barrett and Morton in \cite{BarMorCMAME1984}, see also \cite{DemOdeJCOMPUTPHYS1986,DemOdeCMAME1986,LouHugFraCMAME1987,FraHugLouMirNUMMAT1987,LouFraHugMirCMAME1987}. More precisely, if a fully discontinuous finite element space is used for the test functions, as in the DPG method of Bottasso et al.~\cite{BotMicSacCMAME2002}, the optimal test functions can be approximated locally in an enriched finite element space. The same procedure can also be used to compute local \emph{a posteriori} error estimates to guide automatic adaptive mesh refinement as demonstrated in \cite{Demkowicz2011} in the context of convection-dominated diffusion problems and in \cite{Niemi2011b} in the context of shell boundary layers.

Previous work in the DPG framework applied to Problem~\eqref{eqn:the equation}, has utilized adaptive mesh refinements \cite{DemBOOK2007,Demkowicz2011} to fully resolve fine scales (boundary and interior layers) in the solution. However, in engineering applications this is not always possible. In this work, we focus on the coarse-mesh accuracy of the DPG method. That is, we investigate how the effect of fine scales to the solution can be taken into account without resolving them by the trial space mesh. It should be noted that coarse-mesh accuracy is valuable even in the context of adaptive mesh refinement to avoid unnecessary mesh refinement.

We complement the earlier works \cite{DemGopNUMPDE2010,Demkowicz2011} by resolving the stable test function space with respect to the optimal (quasi-optimal) test space norm. The corresponding DPG method is expected to be uniformly stable with respect to vanishing diffusion, but the local problems for the test functions become singularly perturbed and therefore more difficult to solve. The main contribution of our study is to show how this problem can be approached by utilizing a carefully designed element sub-grid discretization. The present work is a continuation of our recent work \cite{Niemi2011a,Niemi2011} and extends the earlier analysis in one space dimension to two dimensions. 

The paper is structured to be self-contained and focused on bridging the theoretical framework to a practical understanding. Section \ref{sec:abstract framework} lays down the abstract functional analytic setup of the method. Then in Section \ref{sec:application to the equation}, we specialize the theoretical framework to Problem \eqref{eqn:the equation} and introduce different test space norms. In Section \ref{sec:algorithmic considerations}, we re-explain the DPG framework from a purely algorithmic point of view in order to highlight the practical implications of the different choices of the test space norm. We provide details on how to construct the linear algebraic systems constituting the local optimal test functions. Numerical results are presented in Section \ref{sec:numerical results} and the paper ends with conclusions and remarks in Section \ref{sec:conclusions}.

\section{Automatically Stable DPG Framework} \label{sec:abstract framework}
\subsection{Background}
The starting point of the formulation is an abstract variational statement of the form: Find $\vu \in \cU$ such that
\begin{equation} \label{eqn:AVBVP}
\cB(\vw,\vu) = \cL(\vw), \quad \forall \vw \in \cW
\end{equation}
where the test function space $\cW$ and the trial function space $\cU$ are assumed to be Hilbert spaces under inner products $\scp{\cdot}{\cdot}_{\cW}$ and $\scp{\cdot}{\cdot}_{\cU}$ with the corresponding norms $\norm{\cdot}_{\cW}$ and $\norm{\cdot}_{\cU}$. We shall assume that the spaces are real because our application, Problem~\eqref{eqn:the equation}, deals only with real-valued functions. 

If the spaces $\cW$ and $\cU$ are the same and $\cB$ induces an inner product, then the standard Galerkin method delivers the best approximation of the exact solution in the corresponding norm known as the energy norm. In the case of Problem~\eqref{eqn:the equation} the conventional weak formulation is associated with the bilinear and linear forms
\begin{equation} \label{eqn:standard variational form}
\cB(v,u) = \int_{\Omega} \left( \epsilon \vnabla v \cdot \vnabla u + v  (\va \cdot \vnabla u) \right)\,\mathrm{d} \Omega, \quad \cL(v) = \int_{\Omega} vf\,\mathrm{d}\Omega
\end{equation}
Here $\cB$ does not induce an inner product because the second term is not symmetric. Consequently the best approximation property in the energy norm is lost together with numerical stability, see \cite{BroHugCMAME1982,HugScoFraECM2007}. This has motivated the use of the Petrov-Galerkin method, where the combination of specially engineered test functions with added stabilization terms in the variational formulation, has improved the numerical approximations considerably. The DPG method discussed here restores the best approximation property by the computation of a special set of test functions for given trial functions that produce a symmetric positive-definite algebraic system when substituted in the bilinear form. The method can be explored from many different perspectives. It may be viewed as an automatic apparatus designed for guaranteeing the famous inf-sup condition of the Babu\v{s}ka-Brezzi theory of discrete stability, see \cite{BabNUMERMATH1971,BreRAIRO1974}, or it may be interpreted as a generalized least squares method, see \cite{DemGopIACM2011}. We treat the methodology as a numerical construction which generates a true inner-product structure for any well-posed variational problem.

\subsection{Petrov-Galerkin Method with Optimal Test Functions}
The goal is to find a linear operator $\mT$ which maps the trial function space $\cU$ to the test function space $\cW$ such that
\begin{equation} \label{eqn:goal}
\cB(\mT\vv,\vu) = \cA(\vv,\vu) \quad \forall \vv, \vu \in \cU
\end{equation}
where $\cA(\vv,\vu)$ is a symmetric coercive bilinear form on $\cU \times \cU$. If $\{ \ve_1,\ve_2,\ldots,\ve_n\}$ denotes a set of trial functions that spans a subspace $\cU_n$ of $\cU$, then the use of the test functions $\mT \ve_i$ in the Petrov-Galerkin method
\begin{equation} \label{eqn:PG method}
\cB(\mT\ve_i,\vu_n) = \cL(\mT\ve_i), \quad i = 1,\dots,n
\end{equation}
defines $\vu_n$ as the orthogonal projection of the exact solution $\vu$ onto the trial space $\cU_n$ with respect to the inner product $\cA(\cdot,\cdot)$:
\[
\cA(\vv_n,\vu-\vu_n) = 0 \quad \forall \vv_n \in \cU_n
\]
In other words, $\vu_n$ is the best approximation of $\vu$ in $\cU_n$ in the corresponding `energy' norm, that is,
\[
\enorm{\vu - \vu_n}_{\cU} \leq \enorm{\vu - \vv}_{\cU} \quad \forall \vv \in \cU_n
\]
where
\[
\enorm{\vu}_{\cU} = \sqrt{\cA(\vu,\vu)} \quad \forall \vu \in \cU
\]

A canonical way to define the \emph{trial-to-test operator} $\mT:\cU \rightarrow \cW$ is to define $\mT\vu \in \cW$ through the variational equality
\begin{equation} \label{eqn:trial-to-test}
\scp{\vw}{\mT\vu}_{\cW} = \cB(\vw,\vu) \quad \forall\, \vw \in \cW
\end{equation}
so that we have
\[
\cB(\mT\vv, \vu) = \scp{\mT\vv}{\mT\vu}_{\cW}
\]
and thus achieve our goal with 
\begin{equation} \label{eqn:energy norm}
\cA(\vv,\vu) = \scp{\mT\vv}{\mT\vu}_{\cW}, \quad  \enorm{\vu}_{\cU} = \norm{\mT \vu}_{\cW}
\end{equation}

\subsection{Localization}
For each trial function $\ve_i$, $i=1,\ldots,n$, Equation~\eqref{eqn:trial-to-test} constitutes an auxiliary, infinite-dimensional variational problem associated with the test function space $\cW$. The standard Galerkin method may be used to approximate solutions to Equation~\eqref{eqn:trial-to-test}. However, for typical variational formulations, like for the one corresponding to the forms of \eqref{eqn:standard variational form}, Equation~\eqref{eqn:trial-to-test} is from the computational perspective at least as `heavy' as the original problem. 

The new insight of Demkowicz and Gopalakrishnan, see~\cite{DemGopNUMPDE2010},  was to formulate the variational problem \eqref{eqn:AVBVP} in an \emph{ultra-weak hybrid form} in such a way that the test function space $\cW$ becomes fully discontinuous on a partitioning of $\Omega$ into elements $\{K_1,K_2,\ldots,K_N\}$:
\[
\cW = \{ \vw : \vw\!\!\mid_{K} \in \cW(K) \; \text{for every element $K$}\}
\] 
As long as the inner product $\scp{\cdot}{\cdot}_{\cW}$ does not couple test functions defined on individual elements, the auxiliary problems defined by Equation~\eqref{eqn:trial-to-test} can be solved locally at the element level. The localized version of the problem can be written in terms of the broken forms
\[
\scp{\cdot}{\cdot}_{\cW} = \sum_{K} \scp{\cdot}{\cdot}_{\cW(K)}^2, \quad \cB(\cdot,\cdot) = \sum_{K} \cB_K(\cdot,\cdot)
\]
which allows the global problem \eqref{eqn:trial-to-test} to be written as finding $\mT \vu \in \cW(K)$ such that
\begin{equation} \label{eqn:local trial-to-test}
\scp{\vw}{\mT\vu}_{\cW(K)} = \cB_K(\vw,\vu) \quad \forall \vw \in \cW(K)
\end{equation}
This problem may then be solved approximately using the standard Galerkin method in an enriched finite element space $\cW_{\tilde{n}}(K) \subset \cW(K)$, that is, find $\tilde{\mT} \vu \in \cW_{\tilde{n}}(K)$ such that
\begin{equation} \label{eqn:approximate local trial-to-test}
\scp{\vw}{\tilde{\mT}\vu}_{\cW(K)} = \cB_K(\vw,\vu) \quad \forall\vw \in \cW_{\tilde{n}}(K)
\end{equation}

\subsection{Well-posedness and Robustness} \label{sec:robustness}
The above construction can be carried out provided that the original variational problem \eqref{eqn:AVBVP} is well-posed, that is, provided that the bilinear form $\cB(\vw,\vu)$ satisfies the conditions of the Babu\v{s}ka-Necas-Nirenberg Theorem, which are (see~\cite{BabNUMERMATH1971} and \cite{BreRAIRO1974,DemICES2006})
\begin{align}
\label{eqn:continuity}
&\cB(\vw,\vu) \leq C_1 \norm{\vw}_{\cW} \norm{\vu}_{\cU} \quad \forall \vw \in \cW, \; \vu \in \cU \\
\label{eqn:surjectivity of T}
&\sup_{\vu \in \cU} \frac{\cB(\vw,\vu)}{\norm{\vu}_{\cU}} \geq C_2 \norm{\vw}_{\cW}  \\
\label{eqn:surjectivity of T*}
&\sup_{\vw \in \cW} \frac{\cB(\vw,\vu)}{\norm{\vw}_{\cW}} \geq C_3 \norm{\vu}_{\cU} 
\end{align}
The first condition implies that the right hand side of Equation~\eqref{eqn:trial-to-test} is a bounded linear functional in $\cW$ for every $\vu \in \cU$. Therefore $\mT \vu \in \cW$ exists and is unique.

The second condition guarantees that the mapping $\mT$ is surjective, that is, every $\vw \in \cW$ can be expressed in the form $\vw = \mT\vu$ for some $\vu \in \cU$. Namely, if this was not the case, there would exist a non-zero $\tilde{\vw} \in \cW$ such that $\scp{\tilde{\vw}}{\mT\vu}_{\cW} = 0$ for all $\vu \in \cU$. But then also
\[
\cB(\tilde{\vw},\vu) = 0 \quad \forall \vu \in \cU 
\]
by Equation~\eqref{eqn:trial-to-test} so that Equation~\eqref{eqn:surjectivity of T} implies $\tilde{\vw} = \boldsymbol{0}$, which is a contradiction.

Finally, the coercivity of the bilinear form $\cA$ follows from condition \eqref{eqn:surjectivity of T*}. To see this, we  recall that the norm of an element in a Hilbert space $\cH$ can be computed by duality as
\begin{equation} \label{eqn:duality}
\norm{\vu}_{\cH} = \sup_{\vv \in \cH} \frac{\scp{\vv}{\vu}_{\cH}}{\norm{\vv}_{\cH}}
\end{equation}
so that the energy norm can be expressed as
\[
\enorm{\vu}_{\cU} = \sup_{\vw \in \cW} \frac{\scp{\vw}{\mT\vu}_{\cW}}{\norm{\vw}_{\cW}} = \sup_{\vw \in \cW} \frac{\cB(\vw,\vu)}{\norm{\vw}_{\cW}}
\]
Now Equation~\eqref{eqn:surjectivity of T*} implies that
\[
\cA(\vu,\vu) \geq C_3^2 \norm{\vu}_{\cU}^2 \quad \forall\,\vu \in \cU
\]

The characteristics of the energy norm depend on the values of the constants $C_1$ and $C_3$ which in turn typically depend on the physical parameters and the problem geometry. As a result the Petrov-Galerkin method \eqref{eqn:PG method} may not be robust, that is, uniformly accurate with respect to the parameter values of a particular physical problem that is being modeled, see~\cite{BabSurSJNA1992}. However, the parametric dependence can be removed completely from the bilinear form $\cA(\cdot,\cdot)$ by using a special norm for the test function space in \eqref{eqn:trial-to-test}. This norm is referred to as the \emph{optimal test space norm} in the DPG literature, see for instance \cite{ZitMugDemGopParCalJCP2011}, and can be defined as 
\begin{equation} \label{eqn:optimal test space norm}
\enorm{\vw}_{\cW} = \norm{\mT^* \vw}_{\cU}
\end{equation}
where $\mT^*: \cW \rightarrow \cU$ is a linear mapping such that
\[
\scp{\mT^* \vw}{\vu}_{\cU} = \cB(\vw,\vu) \quad \forall \vu \in \cU
\]
Equation~\eqref{eqn:optimal test space norm} defines a proper norm in $\cW$ again under the assumption that the variational problem is well-posed. Namely, the utilization of the duality argument \eqref{eqn:duality}, and conditions \eqref{eqn:continuity}--\eqref{eqn:surjectivity of T} show that
\[
C_2 \norm{\vw}_{\cW} \leq \norm{\mT^*\vw}_{\cU} \leq C_1 \norm{\vw}_{\cW}
\]
that is, $\enorm{\cdot}_{\cW}$ is an equivalent norm to $\norm{\cdot}_{\cW}$.

Finally, condition \eqref{eqn:surjectivity of T*} implies that $\mT^*$ is surjective so that a subsequent application of the duality argument \eqref{eqn:duality} yields
\begin{equation} \label{eqn:norm equivalence}
\enorm{\vu}_{\cU} = \sup_{\vw \in \cW} \frac{\cB(\vw,\vu)}{\enorm{\vw}_{\cW}} = \sup_{\vw \in \cW} \frac{\scp{\mT^* \vw}{\vu}_{\cU}}{\norm{\mT^*\vw}_{\cU}} = \sup_{\vv \in \cU} \frac{\scp{\vv}{\vu}_{\cU}}{\norm{\vv}_{\cU}} = \norm{\vu}_{\cU}
\end{equation}
Thus, the introduction of \eqref{eqn:optimal test space norm} allows us to obtain a problem for each trial function which yields a robust method (convergent in the $\norm{\cdot}_{\cU}$ norm as stated in 
\eqref{eqn:norm equivalence}). This guaranteed best approximation property is what leads us to call the norm induced by Equation \eqref{eqn:optimal test space norm} optimal. As we will observe in the proceeding Section, the ultra-weak approach allows the expression of the optimal test space norm to be deduced directly from the bilinear form. Additionally, the leading terms can be expressed in a closed, localizable form. 

\section{DPG Method for the Convection-Diffusion Problem} \label{sec:application to the equation}
In this section we follow \cite{DemGopNUMPDE2010,DemGopSJNA2011} and derive an explicit DPG variational formulation of the convection-diffusion equation. The ultra-weak DPG variational form of the problem is associated initially to a partitioning of $\Omega$ into elements denoted by $\Omega_h$. We will refer to the collection of all edges in the mesh as $\partial \Omega_h = \cup_K \partial K$. We will use the standard notation where $L_2(\Omega)$ and $\boldsymbol{L}_2(\Omega)$ denote the spaces of square-integrable scalar- and vector-valued functions over $\Omega$, respectively. Moreover, $H^1(\Omega)$ and $\boldsymbol{H}(\mathrm{div},\Omega)$ stand for the subspaces consisting of functions with square-integrable derivatives:
\[
\begin{aligned}
H^1(\Omega) &= \{ v \in L_2(\Omega) \; : \; \vnabla v \in L_2(\Omega) \} \\
\boldsymbol{H}(\mathrm{div},\Omega) &= \{ \mtau \in \boldsymbol{L}_2(\Omega) \; : \; \vnabla \cdot \mtau \in L_2(\Omega) \}
\end{aligned}
\]

\subsection{Derivation of the Hybrid Ultra-Weak Variational Form}
The procedure begins by rewriting the second-order partial differential equation in~\eqref{eqn:the equation} as a pair of first-order equations
\begin{equation} \label{eqn:first-order system}
\begin{aligned}
\epsilon^{-1} \msigma &= -\vnabla u \\
\vnabla \cdot \msigma + \va \cdot \vnabla u &= f
\end{aligned}
\end{equation}
Here we have chosen to define $\msigma$ as the diffusive flux in order to be consistent with the above cited references. Other option would be to include the advective flux $\va u$ in the first equation defining $\msigma$ as in our earlier works \cite{Niemi2011a,Niemi2011}. In our experience, the difference between the two approaches appears to be mostly cosmetic. In this work we will assume that the ambient fluid is incompressible so that $\vnabla \cdot \va = 0$ and that the problem is formulated in a dimensionless form scaled such that $\abs{\va} \sim 1$.

Upon multiplying the first and second equations in \eqref{eqn:first-order system} by a vector-valued test function $\mtau$ and a scalar-valued test function $v$, respectively, and integrating by parts element-wise, we obtain the following bilinear and linear forms corresponding to the abstract variational statement \eqref{eqn:AVBVP}
\begin{equation} \label{eqn:forms}
\begin{aligned}
\cB(\vw,\vu) &= \scp{\epsilon^{-1}\mtau - \vnabla v}{\msigma}_{\Omega_h} - \scp{\vnabla \cdot \mtau + \va \cdot \vnabla v}{u}_{\Omega_h} + \dual{\mtau \cdot \vn}{\hat{u}}_{\partial\Omega_h} + \dual{\hat{\sigma}_n}{v}_{\partial \Omega_h} \\
\cL(\vw) &= \scp{v}{f}_{\Omega_h}
\end{aligned}
\end{equation}
Here $\vn$ denotes the outward unit normal on each element edge and $\hat{\sigma}_n$ corresponds to the (outer) normal component of the total flux $ (\msigma + \va u) \cdot \vn$. Moreover, the test function is denoted by $\vw = (\mtau,v)$ and the trial function by $\vu = (\msigma,u,\hat{u},\hat{\sigma}_n)$. For smooth functions, the notations $\scp{\cdot}{\cdot}_{\Omega_h}$ and $\dual{\cdot}{\cdot}_{\partial \Omega_h}$ stand for the element-wise computed inner products
\[
\scp{f}{g}_{\Omega_h} = \sum_{K \in \Omega_h} \int_K f g \,\mathrm{d}\Omega, \quad \dual{f}{g}_{\partial \Omega_h} = \sum_{K \in \Omega_h} \int_{\partial K} f g\,\mathrm{d}s
\]

The element integrals $\scp{\cdot}{\cdot}_{\Omega_h}$ in Equation~\eqref{eqn:forms} make sense provided that $\msigma, u, \mtau,v$ as well as $\vnabla \cdot \mtau$ and $\vnabla v$ are square integrable over each $K$.  Thus, the test function space is defined as the broken space
\begin{equation} \label{eqn:broken test space}
\cW = \boldsymbol{H}(\mathrm{div},\Omega_h) \times H^1(\Omega_h)
\end{equation}
where
\[
\begin{aligned}
\boldsymbol{H}(\mathrm{div},\Omega_h) &= \{ \mtau \in \boldsymbol{L}_2(\Omega) \; : \; \vnabla \cdot \mtau \in L_2(K) \quad \forall K \in \Omega_h \} \\
H^1(\Omega_h) &= \{ v \in L_2(\Omega) \; : \; \vnabla v \in L_2(K) \quad \forall K \in \Omega_h \}
\end{aligned}
\]

In the general case, the interface action $\dual{\cdot}{\cdot}_{\partial \Omega_h}$ in \eqref{eqn:forms} must be interpreted as an appropriate duality pairing between suitable generalized function spaces defined on $\partial K$. The right setting can be deduced from the theory of trace operators naturally associated to the spaces $\boldsymbol{H}(\mathrm{div},K)$ and $H^1(K)$. These spaces are formally denoted by $H^{-1/2}(\partial K)$ and $H^{1/2}(\partial K)$ for the normal trace $\mtau\!\cdot \vn\mid_{\partial K}$ and the standard trace $v\!\mid_{\partial K}$, respectively.

This leads to the definition of the trial function space as
\begin{equation} \label{eqn:trial space}
\cU = \boldsymbol{L}^2(\Omega) \times L^2(\Omega) \times H_D^{1/2}(\partial \Omega_h) \times H^{-1/2}(\partial \Omega_h)
\end{equation}
where the spaces for the interface variables $\hat{u}$ and $\hat{\sigma}_n$ are defined as
\[
\begin{aligned}
H_D^{1/2}(\partial \Omega_h) &= \{ v\!\mid_{\partial\Omega_h} \; : \; v \in H_D^1(\Omega) \} \\
H^{-1/2}(\partial \Omega_h) &= \{ \veta\!\mid_{\partial\Omega_h} \; : \; \veta \in H(\mathrm{div},\Omega) \}
\end{aligned}
\]
Here $H_D^1(\Omega)$ refers to the subspace of $H^1(\Omega)$ encompassing the Dirichlet boundary condition $u=g$ on $\partial\Omega$. 
The topology of the numerical trace spaces can be characterized for instance by using the extensions with minimal norm:
\[
\begin{aligned}
\norm{\hat{u}}_{H^{1/2}(\partial \Omega_h)} &= \{ \inf \norm{v}_{H^1(\Omega)}  \; : \; v \in H^1(\Omega) \; \text{such that} \; v\!\mid_{\partial \Omega_h} = \hat{u} \} \\
\norm{\hat{\sigma}_n}_{H^{-1/2}(\partial \Omega_h)} &= \{ \inf \norm{\veta}_{\boldsymbol{H}(\mathrm{div},\Omega)} \; : \; \veta \in \boldsymbol{H}(\mathrm{div},\Omega)  \; \text{such that} \; \veta \cdot \vn\!\mid_{\partial \Omega_h} = \hat{\sigma}_n \}
\end{aligned}
\]
More details regarding these definitions can be found from \cite{DemGopSJNA2011}. For our purposes it suffices to note that piecewise polynomial functions in $H_D^{1/2}(\partial\Omega_h)$ must be globally continuous on $\partial \Omega_h$ whereas discontinuities at the vertices of $\Omega_h$ are allowed in the space $H^{-1/2}(\partial \Omega_h)$.

\subsection{Finite Element Trial Spaces}
We refer to $P_t$ as the set of polynomials of degree lower or equal to $t$ and to the corresponding tensor product family by $Q_{p_x,p_y} = P_{p_x} \times P_{p_y}$ and set $Q_p = Q_{p,p}$. The trial space is defined then as
\[
\cU_n = S_h^2 \times S_h \times M_h \times N_h
\]
where
\[
\begin{aligned}
S_h &= \{ v \in L_2(\Omega) : v\!\!\mid_K \; \in Q_p \; \text{for every element $K$} \} \\
M_h &= \{ v \in H_D^1(\Omega) : v\!\!\mid_E \; \in P_{p+1} \; \text{for every edge $E$} \} \\
N_h &= \{ v : v\!\!\mid_E \; \in P_{p} \; \text{for every edge} \; E \}
\end{aligned}
\]
In this work we use the Bernstein polynomials as a basis for $P_t$. The $t+1$ basis polynomials are defined as
\[
B_{i,t}(x) = \binom{t}{i} x^i(1-x)^{t-i}, \quad i = 0,1,\ldots,t
\]

\subsection{Localization in the Convection-Diffusion Problem} \label{sec:localization}
Upon selection of a norm $\norm{\cdot}_{\cW}^2 = \scp{\cdot}{\cdot}_{\cW}$ for the broken test function space $\cW$ defined in~\eqref{eqn:broken test space}, a precise form of the localized auxiliary problem~\eqref{eqn:local trial-to-test} is obtained. This problem is solved approximately using the standard Galerkin method and an enriched finite element space defined on each element $K$. We will associate this space to a sub-partitioning of $K$ into elements~$\{\tilde{K}_1,\ldots,\tilde{K}_{\tilde{N}}\}$.

More precisely, we define an $\boldsymbol{H}(\mathrm{div},K) \times H^1(K)$-conforming finite element space for the test functions as
\[
\cW_{\tilde{n}}(K) = \cT_{\tilde{h}} \times \cV_{\tilde{h}}
\]
where
\[
\begin{aligned}
\cT_{\tilde{h}} &= \{ \mtau \in H(\mathrm{div},K) : \mtau \!\!\mid_{\tilde{K}} \in Q_{\tilde{p}+1,\tilde{p}} \times Q_{\tilde{p},\tilde{p}+1} \; \text{for every sub-element $\tilde{K}$} \} \\
\cV_{\tilde{h}} &= \{ v \in H^1(K) : v \!\!\mid_{\tilde{K}} \in Q_{\tilde{p}} \; \text{for every sub-element $\tilde{K}$} \}
\end{aligned}
\]
Here we denote $\tilde{p} = p + \Delta p$ so that the level of enrichment is controlled by both the number of sub-elements $\tilde{N}$ and the increment in the order of approximation $\Delta p$. In what follows, we will describe three different choices of the test space norm and discuss appropriate enrichments for resolving the corresponding test functions. 

\subsubsection*{The Standard Test Space Norm}
A norm known as the \emph{standard norm}, or mathematician's norm, may be constructed by composing the natural norms of the spaces $\boldsymbol{H}(\mathrm{div}, \Omega_h)$ and $H^1(\Omega_h)$ as in
\[
\norm{\vw}_{\mathrm{SN}}^2 = \norm{\vnabla \cdot \mtau}_{\Omega_h}^2 + \norm{\mtau}_{\Omega_h}^2 
+ \norm{\vnabla v}_{\Omega_h}^2 + \norm{v}_{\Omega_h}^2
\]
This choice leads to two decoupled, standard variational problems for $\mtau$ and $v$ over each element $K$ which can be solved sufficiently accurately by using an enriched finite element space with $\Delta p = 2$ and $\tilde{N} = 1$, see \cite{DemGopNUMPDE2010}.

\subsubsection*{The Weighted Test Space Norm}
While the standard norm provides a local problem that is easy to solve, it leads to approximations of $u$ that underestimate the true values considerably. It was determined in \cite{Demkowicz2011}, that the use of a special weighting of the norm near the inflow boundary could improve the accuracy. The original version of the weighted norm, as introduced in \cite{Demkowicz2011}, is based directly on the standard norm and reads
\[
\norm{\vw}^2_{\mathrm{WN}} = \norm{\vnabla \cdot \mtau}_{\Omega_h,\beta}^2 + \norm{\mtau}_{\Omega_h,\beta}^2 
+ \norm{\vnabla v}_{\Omega_h,\beta}^2 + \norm{v}_{\Omega_h,\beta}^2
\]
where $\beta$ is a piecewise constant norm weight such that
\[
\beta(\vx) = \begin{cases}
\gamma, &\text{if $d(\vx,\Gamma^{-}) \leq \delta$ \& $d(\vx,\Gamma^{+}) \geq \delta$} \\
1, &\text{otherwise}
\end{cases}
\]
where $d(\vx,\Gamma^{-})$ and $d(\vx,\Gamma^{+})$ stand for the (normal) distances between the point $\vx$ and the inflow and outflow/no-flow boundaries, respectively. These are defined as
\[
\begin{aligned}
\Gamma^{-} &= \{ \vx \in \partial\Omega : \va \cdot \vn < 0 \} \\
\Gamma^{+} &= \{ \vx \in \partial\Omega : \va \cdot \vn \geq 0 \}
\end{aligned}
\]
It should be noted that different refined variants of the weighted norm have been recently introduced and studied in \cite{DemHeuICES2011}. These refinements have resemblance to the quasi-optimal norm introduced next.

\subsubsection*{The Quasi-Optimal Test Space Norm}
While the weighted test space norm resolves some of the deficiencies of the standard norm, the approach is not optimal in the sense of the DPG theory as outlined in Section \ref{sec:robustness}. In contrast to the weighted test space norm, we describe here the optimal test space norm \eqref{eqn:optimal test space norm}. However, the true optimal norm leads to a series of global problems so that some numerical modifications are necessary to make the approach suitable for practical computations. This way we arrive to a quasi-optimal test space norm which retains the essential features of the optimal test space norm but is localizable. 

Let the trial space norm be the natural norm arising from the definition \eqref{eqn:trial space}
\[
\norm{\vu}_{\cU} = \left\{\norm{\msigma}_{L_2(\Omega)}^2 + \norm{u}_{L_2(\Omega)}^2 + \norm{\hat{u}}_{H^{1/2}(\partial \Omega_h)}^2 + \norm{\hat{\sigma}_n}_{H^{-1/2}(\partial \Omega_h)}^2\right\}^{1/2}
\]
Thus, using the definition \eqref{eqn:forms} of the bilinear form, we can express the optimal test space norm \eqref{eqn:optimal test space norm} as
\[
\enorm{\vw}_{\cW}^2 = \norm{\epsilon^{-1}\mtau - \vnabla v}_{\Omega_h}^2 + \norm{\vnabla \cdot \mtau + \va \cdot \vnabla v}_{\Omega_h}^2 + \norm{[\mtau \cdot \vn]}_{\partial \Omega_h}^2 + \norm{[v\vn]}_{\partial \Omega_h}^2
\]
where
\[
\norm{[\mtau \cdot \vn]}_{\partial \Omega_h} = \sup_{\hat{u} \in H_D^{1/2}(\partial \Omega_h)} \frac{\dual{\hat{u}}{\mtau \cdot \vn}_{\partial \Omega_h}}{\norm{\hat{u}}_{H^{1/2}(\partial \Omega_h)}}, \quad
\norm{[v \vn]}_{\partial \Omega_h} = \sup_{\hat{\sigma}_n \in H^{-1/2}(\partial \Omega_h)} \frac{\dual{v}{\hat{\sigma}_n}_{\partial \Omega_h}}{\norm{\hat{\sigma}_n}_{H^{-1/2}(\partial \Omega_h)}}
\]
The last two terms can be made more precise, see \cite{DemGopSJNA2011}, but the important point is that they involve the `jumps' of $\mtau$ and $v$ along the inter-element boundaries and therefore prevent the element-wise computation of the optimal test functions. However, a localizable variant of the optimal norm may obtained simply by dropping these terms and adding an integral contribution of $v$, so as to remove the `zero-energy mode' $(\mtau,v) = (\boldsymbol{0},\mathrm{const.})$. This leads us to the so-called \emph{quasi-optimal test space norm}
\begin{equation} \label{eqn:quasi-optimal test space norm}
\norm{\vw}_{\mathrm{QON}}^2 = \norm{\epsilon^{-1}\mtau - \vnabla v}_{\Omega_h}^2 + \norm{\vnabla \cdot \mtau + \va \cdot \vnabla v}_{\Omega_h}^2 + \alpha_1 \norm{\mtau}_{\Omega_h}^2 + \alpha_2 \norm{v}_{\Omega_h}^2 
\end{equation}
where $\alpha_1 \geq 0$ and $\alpha_2 > 0$ are numerical regularization parameters to be chosen. Notice that while the addition of a contribution of $\norm{\mtau}_{\Omega_h}^2$ is not necessary, we have observed that a proper regularization of $\mtau$ improves the accuracy of the final approximation. In our algorithm we choose $\alpha_1 = \epsilon^{-3/2}$ and  $\alpha_2 = 1$.

A critical issue in the application of the quasi-optimal norm \eqref{eqn:quasi-optimal test space norm} is that the auxiliary problems for the test functions become singularly perturbed. This was explored in the one-dimensional case in \cite{Niemi2011a,Niemi2011} where a specially designed sub-grid mesh was employed to resolve the potential boundary layers in the test functions. To extend the idea to two spatial dimensions, we observe that the strong form of the problem \eqref{eqn:local trial-to-test} associated to the quasi-optimal test space norm \eqref{eqn:quasi-optimal test space norm} takes the form
\begin{equation} \label{eqn:optimal test function problem}
\left\{
\begin{aligned}
-\vnabla (\vnabla \cdot \mtau) + (\epsilon^{-2} + \alpha_1) \mtau - \vnabla (\va \cdot \vnabla v) - \epsilon^{-1} \vnabla v &= \varepsilon^{-1}\msigma + \vnabla u \\
- \Delta v - \va \cdot \vnabla (\va \cdot \vnabla v) + \alpha_2 v -\va \cdot \vnabla (\vnabla \cdot \mtau) + \epsilon^{-1} \vnabla \cdot \mtau &= \vnabla \cdot \msigma + \va \cdot \vnabla u
\end{aligned}
\right.
\end{equation}
in every $K$. Notice that the system \eqref{eqn:optimal test function problem} corresponds to the Euler-Lagrange equations for the variational problem \eqref{eqn:trial-to-test} and is accompanied with inhomogeneous natural boundary conditions involving also the interface variables $\hat{\sigma}_n$ and $\hat{u}$.

\begin{figure}
\begin{center}
\includegraphics{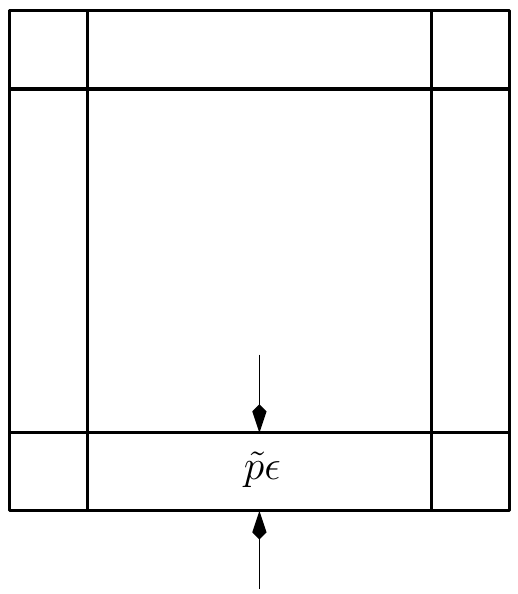}
\end{center}
\caption{A $3 \times 3$ Shishkin mesh to resolve the boundary layers of the optimal test functions.}
\label{fig:subgrid}
\end{figure}

While a detailed regularity and asymptotic analysis of the system \eqref{eqn:optimal test function problem} is out of the scope of this paper, a preliminary examination of the system when $K = (0,1) \times (0,1)$ and $\va = (1,0)$ or $\va = (0,1)$ reveals that the solution may exhibit regular, exponential boundary layers behaving like
\[
e^{-\lambda x_1}, \quad e^{-\lambda (1-x_1)}, \quad e^{-\lambda x_2}, \quad e^{-\lambda (1-x_2)}
\]
where the characteristic exponent is given by $\lambda = \epsilon^{-1}\sqrt{1 + \alpha_1\epsilon^2}$. Notice that when $\alpha_1 = \epsilon^{3/2}$, we have $\lambda = \epsilon^{-1} + \cO(\epsilon^{-1/2})$. To account for these features in the automatic computation of the optimal test functions, we employ a layer-adapted Shishkin mesh on each $K$ as shown in Figure \ref{fig:subgrid}. Such meshes are commonly used in the approximation of singularly perturbed problems, see \cite{LinBOOK2010,ShiRJNAMM1989}. We follow a strategy which consists of adding one layer of 'needle' elements of width $\cO(\tilde{p} \epsilon)$ near the boundary and has been proven to be optimal in the context of the $hp$-version of FEM for reaction-diffusion equations, see \cite{SchSurMATHCOMP1996,XenM3AS1998,MelSchSJNA1998}.

We emphasize that the sub-grid discretization is needed only when using the quasi-optimal test space norm. Failure to address the fine-scale features of the corresponding test functions (here using the sub-grid) degrades the accuracy of the final solution and dissipates all benefits of the quasi-optimal test norm.

%

\section{Algorithmic Considerations} \label{sec:algorithmic considerations}
%
At the high level, the DPG method with optimal test functions is similar to the standard finite element method detailed in Algorithm~\ref{a:dpg}. Initially, the mesh information is input and mappings from the local stiffness matrix to the global matrix are determined. After this, we loop over the elements and form the element contributions to the
global stiffness matrix and load vector. The main conceptual difference between the DPG method and standard FEM is the following: the local stiffness matrix will be computed as a by-product of computing the optimal test functions. Another difference is that, because of the hybridization, the trial spaces are non-standard and require some specialized handling. The remaining portion of this section will describe both differences in more detail. For a review of standard FEM implementation, see \cite{Hughes1987,DemBOOK2007}.
\begin{algorithm}
\caption{The DPG finite element method with optimal test functions}\label{a:dpg}
\begin{algorithmic}[1]
  \State Initialize mesh
  \State Compute all local to global maps
  \ForAll{elements in mesh}
    \State Compute local contributions $\mathbf{K}_\ell,\mathbf{F}_\ell$ using \textsc{OptimalTestFunctions} (Algorithm~\ref{a:magic})
    \State Assemble $\mathbf{K}_\ell\rightarrow \mathbf{K}$ and $\mathbf{F}_\ell\rightarrow \mathbf{F}$
    \State Apply boundary conditions
  \EndFor
  \State Solve $\mathbf{K}\mathbf{u}=\mathbf{F}$
\end{algorithmic}
\end{algorithm}
\subsection{Selection of the Trial Spaces}
The space of trial functions in the DPG method consists of basis functions used to discretize variables on the element
interiors ($S_h$ in our application), as well on the edges ($M_h,N_h$). Notice that on edges we must span spaces where corner degrees of freedom are assembled with those of neighboring edges ($M_h$) as well as
spaces where corner degrees of freedom are independent ($N_h$). 

What this means from the implementation standpoint, is that degrees of freedom from each space must be given unique numbers which correspond to their location in the global stiffness matrix. This is not straightforward because the basis functions are defined over varying topological objects (element interiors/edges) and discretize different variables, a consequence of using a first-order system. A flexible finite element framework is needed to address these issues. It is the authors' experience that even in the context of rectangular domains and structured grid meshes that this is a non-trivial task which is not immediately appreciated when digesting the DPG framework.

\subsection{Computation of the Optimal Test Functions}
The optimal test function problem~\eqref{eqn:approximate local trial-to-test} is solved locally at the element level. By specifying a trial function $\ve_i \in \cU_n$ one can compute the corresponding optimal test function $\vw = \tilde{\mT}\ve_i\!\!\mid_K \in \cW_{\tilde{n}}(K)$ using the standard Galerkin method to solve
\[
\scp{\delta \vw}{\vw}_{\cW(K)} = \cB_K(\delta \vw,\ve_i) \quad \forall \delta \vw \in \cW_{\tilde{n}}(K) \subset \cW(K)
\]
Notice that since the trial function $\ve_i$ is given, the bilinear form $\cB_K(\delta \vw,\ve_i)$ induces a linear functional in $\cW(K)$ and is used to form the right-hand side of the optimal test function problem. Notice that each trial function will lead to a different right-hand side in the above problem. We will denote by $\mF_{\mathrm{opt}}$ the matrix whose columns are the load vectors corresponding to each trial function. The bilinear form on the left hand side is defined by the selected test space inner product and is independent of the trial functions. We will denote by $\mK_{\mathrm{opt}}$ the corresponding stiffness matrix. The columns of the solution $\mU_{\mathrm{opt}}$ to the system $\mK_{\mathrm{opt}} \mU_{\mathrm{opt}} = \mF_{\mathrm{opt}}$ represent the approximated optimal test functions.

Algorithm~\ref{a:magic} describes the procedure used to form and solve this linear system. The algorithm is a function of the trial space $\cU_n$,  the enriched finite element space for the optimal test functions $\cW_{\tilde{n}}(K)$, the chosen optimal test function inner product $\scp{\cdot}{\cdot}_{\cW(K)}$, and the bilinear and linear forms $\cB_K(\cdot,\cdot)$ and $\cL(\cdot)$ arising from the problem setup. At this point, we remind the reader that the optimal test functions are vector-valued. When forming the matrix  $\mK_{\mathrm{opt}}$, the choice of optimal test space norm will produce a block matrix structure depending on which  components of optimal test functions interact with each other. Recall that in our application we have three components $\delta\vw=(\delta\tau_1,\delta\tau_2,\delta v)$ and $\vw=(\tau_1,\tau_2,v)$. We will denote the corresponding blocks of the matrix $\mK_{\mathrm{opt}}$ by $[\mK_{\mathrm{opt}}]_{ij}$, $i,j=1,2,3$.

The particular form of the matrix depends on the manner in which the inner product is defined. If we choose the standard norm (SN), the blocks are defined as 
\[
\begin{aligned}
\left[ \mK_{\mathrm{opt}} \right]_{11} &= \scp{\delta \tau_{1,1}}{\tau_{1,1}}_K + (\delta \tau_{1}, \tau_{1})_{K} \\
\left[ \mK_{\mathrm{opt}} \right]_{12} &= (\delta \tau_{1,1}, \tau_{2,2} )_{K} \\ 
\left[ \mK_{\mathrm{opt}} \right]_{13} &= 0 \\
\left[ \mK_{\mathrm{opt}} \right]_{21} &= (\delta \tau_{2,2},\tau_{1,1} )_{K} \\
\left[ \mK_{\mathrm{opt}} \right]_{22} &= (\delta \tau_{2,2}, \tau_{2,2})_{K} + (\delta \tau_{2}, \tau_{2} )_{K} \\
\left[ \mK_{\mathrm{opt}} \right]_{23} &= 0 \\
\left[ \mK_{\mathrm{opt}} \right]_{31} &= 0 \\
\left[ \mK_{\mathrm{opt}} \right]_{32} &= 0 \\
\left[ \mK_{\mathrm{opt}} \right]_{33} &= (\delta v_{,1}, v_{,1})_{K} + \scp{\delta v_{,2}}{v_{,2}}_{K} + (\delta v,v )_{K}
\end{aligned}
\]
where $\scp{\cdot}{\cdot}_K$ denotes the $L_2$ inner product computed over an element $K$ and subscripted comma indicates differentiation with respect to the coordinate indicated by the subscript. If the weighted norm (WN) is used, the structure of the matrix is identical but the $L_2$ inner product is replaced by the appropriately weighted inner product. However, if we choose to work with the quasi-optimal norm (QON), the matrix $\mK_{\mathrm{opt}}$ becomes full:
\[
\begin{aligned}
\left[ \mK_{\mathrm{opt}}\right]_{11} &= (\delta \tau_{1,1},\tau_{1,1})_{K} + (\epsilon^{-2}+\alpha_1) (\delta \tau_{1}, \tau_1)_{K} \\
\left[ \mK_{\mathrm{opt}}\right]_{12} &= (\delta \tau_{1,1},\tau_{2,2} )_{K} \\ 
\left[ \mK_{\mathrm{opt}}\right]_{13} &= -\epsilon^{-1} \scp{\delta \tau_1}{v_{,1}}_K + \scp{\delta\tau_{1,1}}{a_1 v_{,1} + a_2 v_{,2}}_K \\
\left[ \mK_{\mathrm{opt}}\right]_{21} &= (\delta \tau_{2,2}, \tau_{1,1} )_{K} \\
\left[ \mK_{\mathrm{opt}}\right]_{22} &= (\delta \tau_{2,2},\tau_{2,2})_{K} + (\epsilon^{-2}+\alpha_1) (\delta \tau_{2}, \tau_2)_{K} \\
\left[ \mK_{\mathrm{opt}}\right]_{23} &= -\epsilon^{-1} \scp{\delta \tau_2}{v_{,2}}_K + \scp{\delta\tau_{2,2}}{a_1 v_{,1} + a_2 v_{,2}}_K \\
\left[ \mK_{\mathrm{opt}}\right]_{31} &= -\epsilon^{-1} \scp{\delta v_{,1}}{\tau_{1}}_K + \scp{a_1 \delta v_{,1} + a_2 \delta v_{,2}}{\tau_{1,1}}_K \\
\left[ \mK_{\mathrm{opt}}\right]_{32} &= -\epsilon^{-1} \scp{\delta v_{,2}}{\tau_{2}}_K + \scp{a_1 \delta v_{,1} + a_2 \delta v_{,2}}{\tau_{2,2}}_K \\
\left[ \mK_{\mathrm{opt}}\right]_{33} &= (1+a_1^2)\scp{\delta v_{,1}}{v_{,1}}_K + \scp{a_1a_2\delta v_{,1}}{v_{,2}}_K \\
&+ (1+a_2^2)\scp{\delta v_{,2}}{v_{,2}}_K + \scp{a_1a_2\delta v_{,2}}{v_{,1}}_K \\
&+ \alpha_2\scp{\delta v}{v}_K
\end{aligned}
\]

The rows of the right hand side $\mF_{\mathrm{opt}}$ follow the same block structure as $\mK_{\mathrm{opt}}$. The columns correspond to different trial function components in $\vu = (\sigma_1,\sigma_2,u,\hat{u},\hat{\sigma}_n)$. We will denote the corresponding blocks of $\mF_{\mathrm{opt}}$ by $[\mF_{\mathrm{opt}}]_{ij}$, $i=1,2,3$, $j=1,\dots,5$. The actual form of the matrix is determined by the bilinear form $\cB(\cdot,\cdot)$:
\[
\mF_{\mathrm{opt}} =
\begin{bmatrix}
\scp{\epsilon^{-1} \delta \tau_1}{\sigma_1}_{K} &
0 &
-\scp{\delta \tau_{1,1}}{u}_{K} &
\dual{\delta \tau_1 n_1}{\hat{u}}_{\partial K} &
0 \\

0 &
\scp{\epsilon^{-1}\delta \tau_2}{\sigma_2}_{K} &
-\scp{\delta \tau_{2,2}}{u}_{K} &
\dual{\delta \tau_2 n_2}{\hat{u}}_{\partial K} &
0 \\

-\scp{\delta v_{,1}}{\sigma_1}_{K} &
-\scp{\delta v_{,2}}{\sigma_2}_{K} &
-\scp{a_1 \delta v_{,1} + a_2 \delta v_{,2}}{u}_{K} &
0 &
\dual{\delta v}{\hat{\sigma}_n}_{\partial K}
\end{bmatrix}
\]
where $\dual{\cdot}{\cdot}_{\partial K}$ denotes the standard $L_2$ inner product computed over the element boundary $\partial K$ and $n_1,n_2$ denote the components of the unit outward normal to $\partial K$.

The blocks in the matrices $\mK_{\mathrm{opt}}$ and $\mF_{\mathrm{opt}}$ are formed in the standard finite element fashion. Since the matrix $\mathbf{K}_{\mathrm{opt}}$ is mostly dense, we suggest the use of a direct solver with
support for multiple right-hand sides, such as LAPACK's \verb+DGESV+.

Standard wisdom would then suggest that one could compute the local contributions $\mK_{\ell}$ in the natural way, by first looping over the space of optimal test functions and subsequently the corresponding trial functions as in the Petrov-Galerkin method. However, the local matrix
contributions may be computed directly by a matrix-matrix product of
$\mathbf{K}_{\ell}=\mathbf{U}_{\mathrm{opt}}^{T}\mathbf{F}_{\mathrm{opt}}$. This can been seen most easily by  writing the approximate optimal test function as
\[
\tilde{\mT} \ve_i = \sum_{k=1}^{\tilde{n}} [\mU_{\mathrm{opt}}]_{ki} \tilde{\ve}_k
\]
where $\{\tilde{\ve}_1,\dots,\tilde{\ve}_{\tilde{n}}\}$ stands for the basis of the enriched finite element space $\cW_{\tilde{n}}(K)$. Then it follows that
\[
[\mK_{\ell}]_{ij} = \cB(\tilde{\mT} \ve_i, \ve_j) = \sum_{k=1}^{\tilde{n}} [\mU_{\mathrm{opt}}]_{ki} \cB(\tilde{\ve}_k,\ve_j) =  \sum_{k=1}^{\tilde{n}} [\mU_{\mathrm{opt}}]_{ki} [\mF_{\mathrm{opt}}]_{kj}
\]
so that
\[
[\mK_{\ell}]_{ij} = [\mU_{\mathrm{opt}}^T\mF_{\mathrm{opt}}]_{ij}
\]

The local
contribution to the problem right-hand side $\mathbf{F}_l$, however,
must be assembled in the traditional way. From the implementation point of
view, it is advantageous to think of DPG as a routine which computes a
special local matrix which then must be globally assembled.
\begin{algorithm}
\caption{Computation of the optimal test functions}\label{a:magic}
\begin{algorithmic}[1]
\Function{OptimalTestFunctions}{$\cU_n,\cW_{\tilde{n}}(K),\scp{\cdot}{\cdot}_{\cW(K)},\cB_K(\cdot,\cdot),\cL(\cdot)$}
  \ForAll{components $\delta w_i$ of $\delta \vw$}
    \ForAll{components $w_j$ of $\vw$}
      \If{$\delta w_i$ interacts with $w_j$ as defined by $\scp{\cdot}{\cdot}_{\cW(K)}$}
        \State Assemble $[\mathbf{K}_{\mathrm{opt}}]_{ij}$ block of $\mathbf{K}_{\mathrm{opt}}$
      \EndIf
    \EndFor
    \ForAll{components $u_j$ of $\vu$}
      \If{$\delta w_i$ interacts with $u_j$ as defined by $\cB_K(\cdot,\cdot)$}
        \State Assemble $[\mathbf{F}_{\mathrm{opt}}]_{ij}$ block of $\mathbf{F}_{\mathrm{opt}}$
      \EndIf
    \EndFor
  \EndFor
  \State Solve $\mathbf{K}_{\mathrm{opt}}\mathbf{U}_{\mathrm{opt}}=\mathbf{F}_{\mathrm{opt}}$
  \State $\mathbf{K}_\ell = \mathbf{U}_{\mathrm{opt}}^{T}\mathbf{F}_{\mathrm{opt}}$
  \State Assemble $\mathbf{F}_\ell$ by $\cL(\cdot)$
  \State Return $\mathbf{K}_\ell,\mathbf{F}_\ell$
\EndFunction
\end{algorithmic}
\end{algorithm}



\section{Numerical Results} \label{sec:numerical results}
\providecommand{\e}[1]{\ensuremath{\times 10^{#1}}}

In this section we present numerical results which compare the performance of each of the optimal test space norms. While the DPG solution consists of the field variables $\msigma,u$ and the element interface variables $\hat{\sigma}_n, \hat{u}$, our visualizations will show only the $u$ component of the total solution. In our experience it is helpful to visualize all solution variables, particularly while debugging, but the extra plots do not contribute here to a better understanding of the performance of the method. However, we will also show numerical convergence results in terms of the $L_2$ norm for $u$ and $\msigma$ which includes contributions from all solution variables.
%
%

In all of our model problems the computational domain is taken to be
the unit square $\Omega = (0,1) \times (0,1)$. The problems feature
different Dirichlet boundary conditions and different orientations of
the advection vector $\va$ with respect to the mesh.  For each
problem, we compare the coarse-mesh accuracy of the scheme for
different values of the diffusivity ($\epsilon = 10^{-i}$, $i=2,4,6$)
and different test space norm. We shall use the abbreviation SN for
the Standard Norm, WN for the Weighted Norm and QON for the
Quasi-Optimal Norm, see Section \ref{sec:localization} to recall their definitions.

In the test problems we set $p=1$ corresponding to discontinuous,
piecewise bilinear representation of the field variables $\msigma$ and
$u$ on the domain $\Omega$, discontinuous piecewise linear representation of the numerical
flux $\hat{\sigma}_n$ and continuous, piecewise quadratic
representation of the numerical trace $\hat{u}$ on the `skeleton' $\partial \Omega_h$. The trial spaces are
constructed on meshes of ten by ten elements. The optimal test
functions are approximated using a polynomial enrichment value of
$\Delta p=2$.
\subsubsection*{The Eriksson-Johnson problem}
\begin{figure}
\centering
\includegraphics[scale=1.0]{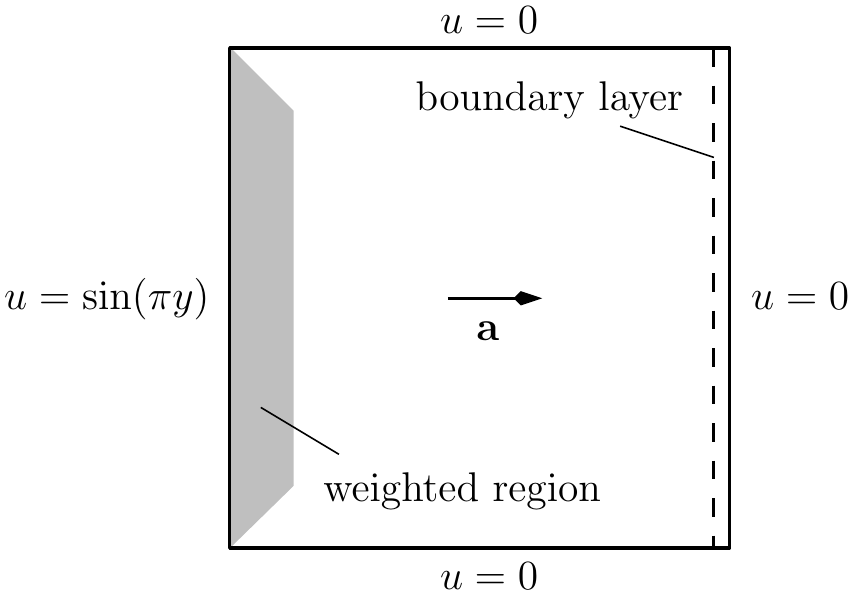}
\caption{Schematics of the Eriksson-Johnson problem. Approximate location of the boundary layer shown as a dashed line and the weighted region is shaded.}
\label{fig:analytic}
\end{figure}

Eriksson and Johnson consider in \cite{EriJohMATHCOMP1993} the equation
\[
u_{,1} - \epsilon \Delta u = f \quad \text{in $\Omega$}\\
\]
representing the convection-diffusion equation with the flow velocity $\va
= (1,0)$. The problem can be approached analytically by separation of
variables. Assuming homogeneous Dirichlet boundary conditions at $x=1$, $y=0$, and $y=1$ and the inflow data $u(0,y) = \sin (\pi y)$, the solution
becomes 
\[
u(x,y) = \frac{\exp\left(\dfrac{(1-s)x}{2\epsilon}\right) - \exp\left(\dfrac{(1+s)x - 2s}{2\epsilon}\right)}{1-\exp\left(-\dfrac{s}{\epsilon}\right)}\sin(\pi y), \quad s = \sqrt{1+  4\pi^2 \epsilon^2}
\]
and features a steep boundary layer of width $\epsilon$ near the
outflow boundary at $x=1$, see Figure~\ref{fig:analytic} for a schematics of the problem. As in all
problem schematics, the dashed line will represent the
location of the boundary layer and the shaded region
represents the area weighted when using the weighted test space norm.

Table \ref{t:analytic} organizes a comparison of the norms for
problems of increasing difficulty in terms of $\epsilon$ versus the element size. A reference solution is
included as a final row which in this case is the analytic solution. From
the first row, one can observe the deficiency of the standard norm as
applied to this problem. Significant undershooting of the exact
solution is observed. This is important even in the context of
adaptivity since error in the coarse scale solutions will lead to
unwanted refinement. The success of the weighted norm in the second
row also becomes clear. The rough character of the true solution is
obtained which would drive adaptivity in the proper locations. However, the
quasi-optimal norm is clearly superior in this case, capturing the
fine scale features on the coarse mesh.
\newcommand{\minmax}[2]{\footnotesize{min/max = #1/#2}}
\def \pathtofig {}
\def \pad {0.225}
\begin{figure}
\caption{Color map plots of $u$ in the Eriksson-Johnson problem with an analytic solution. Color range is uniformly set to $[0,1]$ and the minimum and maximum values of the solution are shown underneath each plot.}\label{t:analytic}
\begin{center}
\begin{tabular}{cccc}
\hline
& $\epsilon = 10^{-2}$ 
& $\epsilon = 10^{-4}$ 
& $\epsilon = 10^{-6}$ 
\\
\hline
\\
\footnotesize{SN} &
\parbox[c]{\pad\textwidth}{\includegraphics[width=\pad\textwidth]{\pathtofig 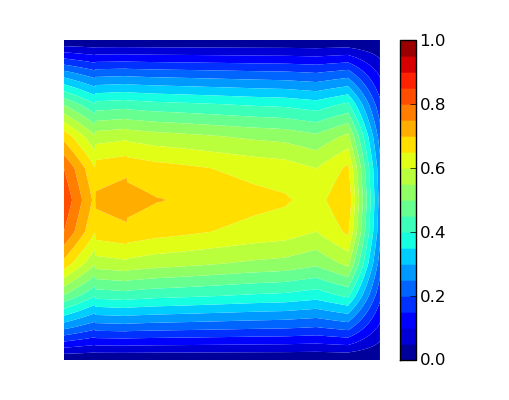}} &
\parbox[c]{\pad\textwidth}{\includegraphics[width=\pad\textwidth]{\pathtofig 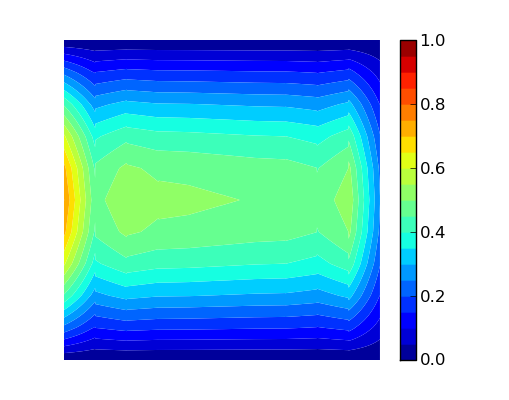}} &
\parbox[c]{\pad\textwidth}{\includegraphics[width=\pad\textwidth]{\pathtofig 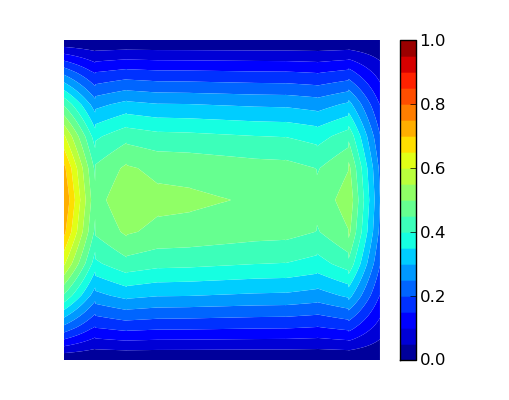}} \\
 & \minmax{0.00}{0.84} & \minmax{0.00}{0.75} & \minmax{0.00}{0.75} \\
\footnotesize{WN} &
\parbox[c]{\pad\textwidth}{\includegraphics[width=\pad\textwidth]{\pathtofig 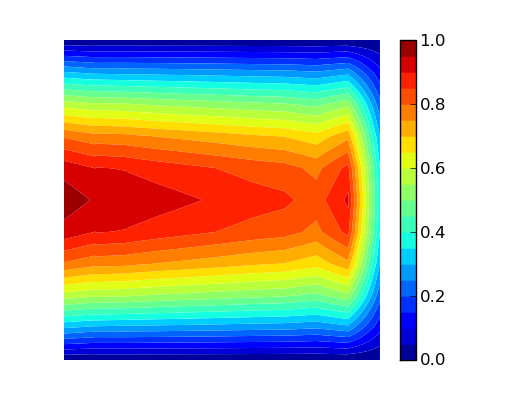}} &
\parbox[c]{\pad\textwidth}{\includegraphics[width=\pad\textwidth]{\pathtofig 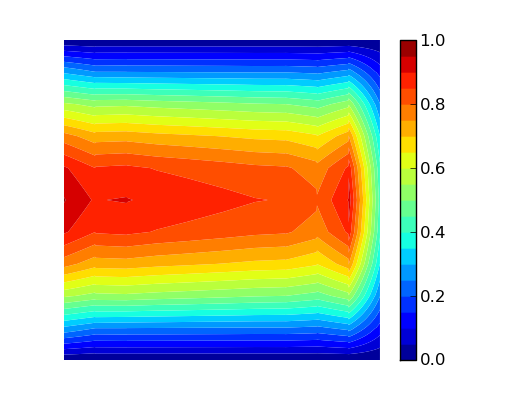}} &
\parbox[c]{\pad\textwidth}{\includegraphics[width=\pad\textwidth]{\pathtofig 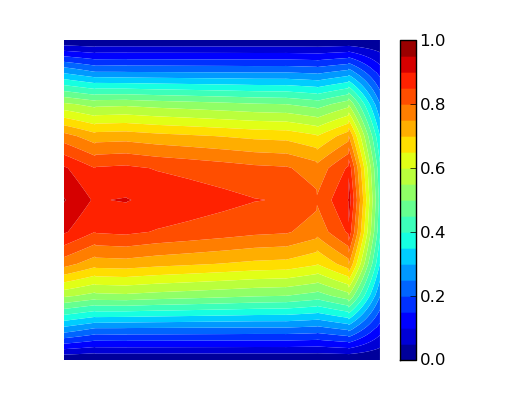}} \\
 & \minmax{0.00}{0.98} & \minmax{0.00}{0.95} & \minmax{0.00}{0.95} \\
\footnotesize{QON} &
\parbox[c]{\pad\textwidth}{\includegraphics[width=\pad\textwidth]{\pathtofig 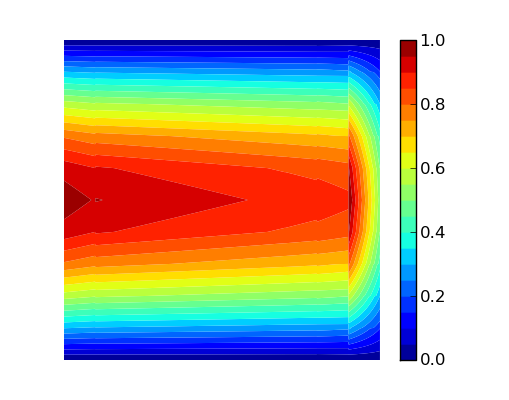}} &
\parbox[c]{\pad\textwidth}{\includegraphics[width=\pad\textwidth]{\pathtofig 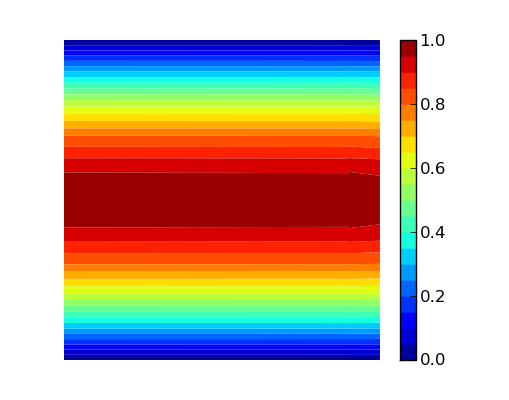}} &
\parbox[c]{\pad\textwidth}{\includegraphics[width=\pad\textwidth]{\pathtofig 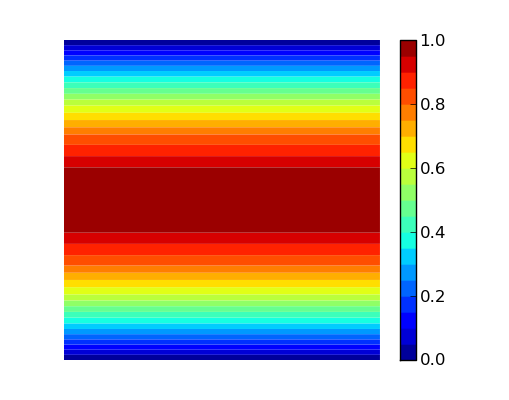}} \\
 & \minmax{0.00}{1.01} & \minmax{0.00}{0.99} & \minmax{0.00}{1.00} \\
\footnotesize{Reference} &
\parbox[c]{\pad\textwidth}{\includegraphics[width=\pad\textwidth]{\pathtofig 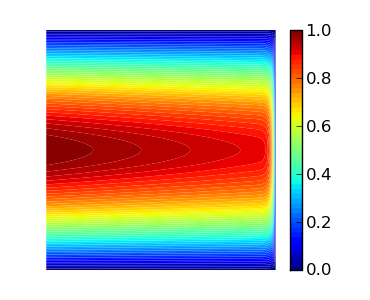}} &
\parbox[c]{\pad\textwidth}{\includegraphics[width=\pad\textwidth]{\pathtofig 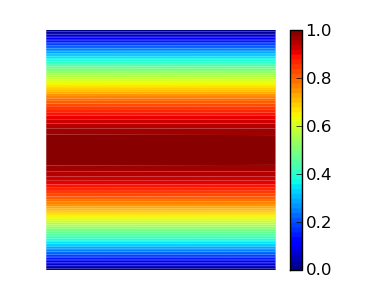}} &
\parbox[c]{\pad\textwidth}{\includegraphics[width=\pad\textwidth]{\pathtofig 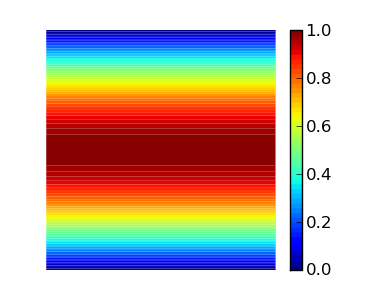}} \\
\\
\hline
\end{tabular}
\end{center}
\end{figure}
\subsubsection*{Advection Skew to the Mesh: Continuous Dirichlet Data}
\begin{figure}
\centering
\includegraphics[scale=1.0]{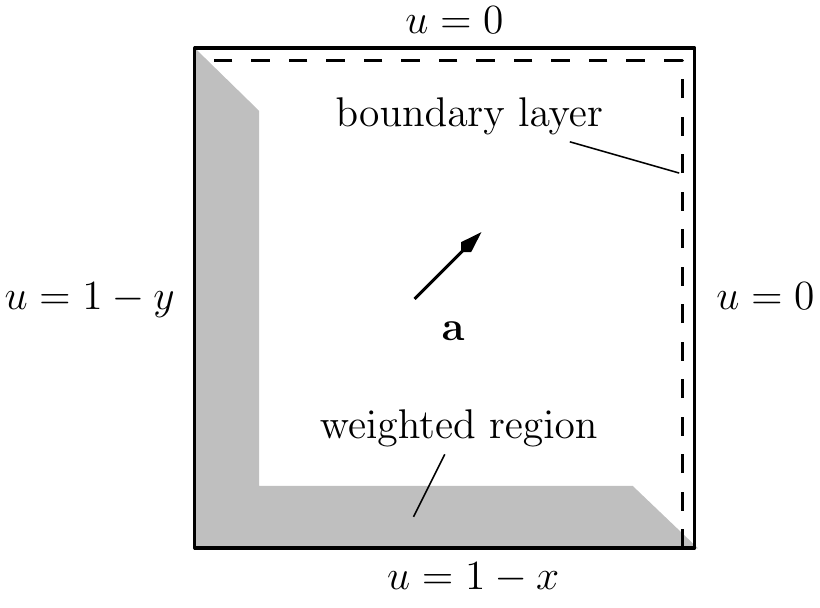}
\caption{Schematics of the advection skew to the mesh problem with continuous Dirichlet data. Approximate location of the boundary layer shown as a dashed line and the weighted region is shaded.}
\label{fig:legal}
\end{figure}

This problem has been studied in
\cite{DemGopNUMPDE2010,Demkowicz2011} and consists of solving
Equation~\eqref{eqn:the equation} with the boundary conditions
\[
u(0,y) = 1-y, \quad u(x,0) = 1-x, \quad
u(x,1) = u(1,y) = 0
\]
and the advection vector $\va = (\cos \theta, \sin \theta)$, $\theta
\in [0^\circ,90^\circ]$. The solution features a boundary layer of
width $\epsilon$ near the outflow boundary. We include a schematic of
the problem in Figure~\ref{fig:legal} where the boundary layer is
indicated by a dashed line and the weighted regions for the WN are
shaded. For this problem, the reference
solutions were computed using standard Galerkin FEM on a
fine mesh which fully resolves the boundary layer. We were not able to
compute such a solution for the problem with $\epsilon=10^{-6}$
although one can compare with the solution for $\epsilon=10^{-4}$ as
the difference in this case is visually indistinguishable. 

\begin{figure}
\caption{Color map plots of $u$ in the advection skew to the mesh problem with continuous Dirichlet data. Color range is uniformly set to $[0,1]$ and the minimum and maximum values of the solution are shown underneath each plot.}\label{t:legal}
\begin{center}
\begin{tabular}{cccc}
\hline
& $\epsilon = 10^{-2}$ 
& $\epsilon = 10^{-4}$ 
& $\epsilon = 10^{-6}$ 
\\
\hline
\\
\footnotesize{SN} &
\parbox[c]{\pad\textwidth}{\includegraphics[width=\pad\textwidth]{\pathtofig 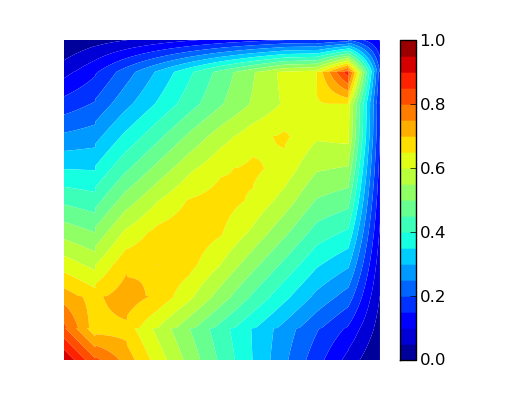}} &
\parbox[c]{\pad\textwidth}{\includegraphics[width=\pad\textwidth]{\pathtofig 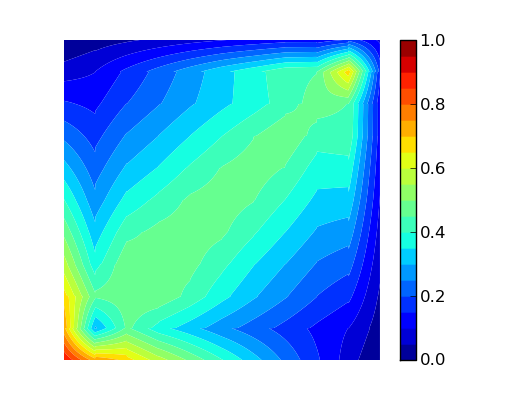}} &
\parbox[c]{\pad\textwidth}{\includegraphics[width=\pad\textwidth]{\pathtofig 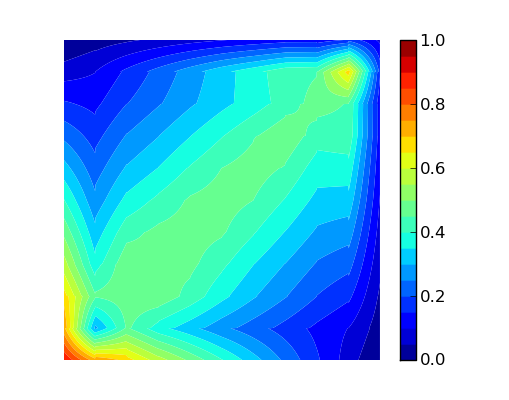}} \\
 & \minmax{0.00}{0.95} & \minmax{0.00}{0.90} & \minmax{0.00}{0.90} \\
\footnotesize{WN} &
\parbox[c]{\pad\textwidth}{\includegraphics[width=\pad\textwidth]{\pathtofig 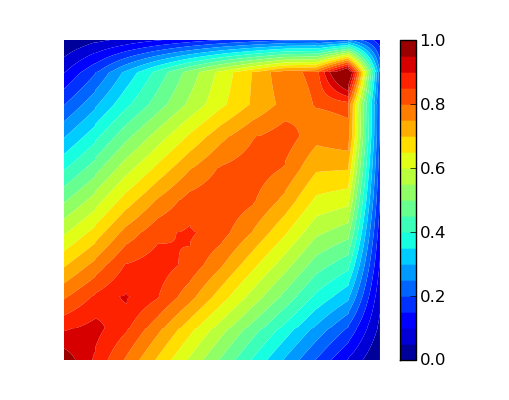}} &
\parbox[c]{\pad\textwidth}{\includegraphics[width=\pad\textwidth]{\pathtofig 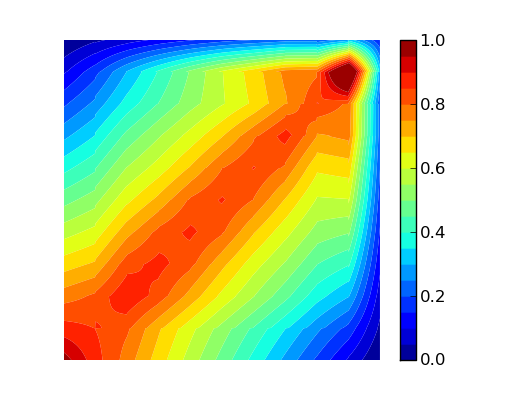}} &
\parbox[c]{\pad\textwidth}{\includegraphics[width=\pad\textwidth]{\pathtofig 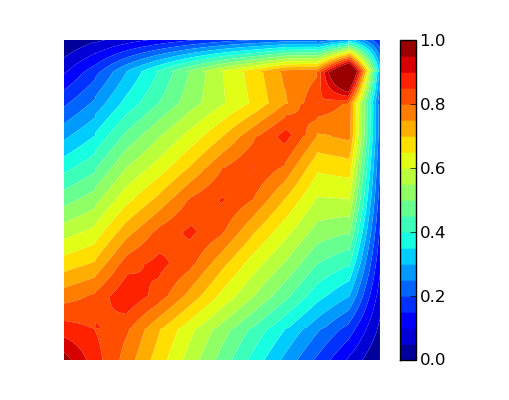}} \\
 & \minmax{0.00}{1.12} & \minmax{0.00}{1.20} & \minmax{0.0}{1.20} \\
\footnotesize{QON} &
\parbox[c]{\pad\textwidth}{\includegraphics[width=\pad\textwidth]{\pathtofig 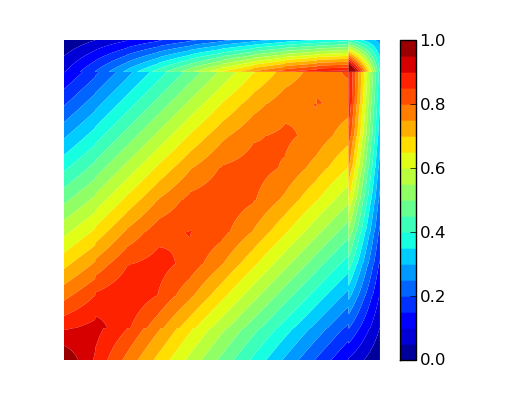}} &
\parbox[c]{\pad\textwidth}{\includegraphics[width=\pad\textwidth]{\pathtofig 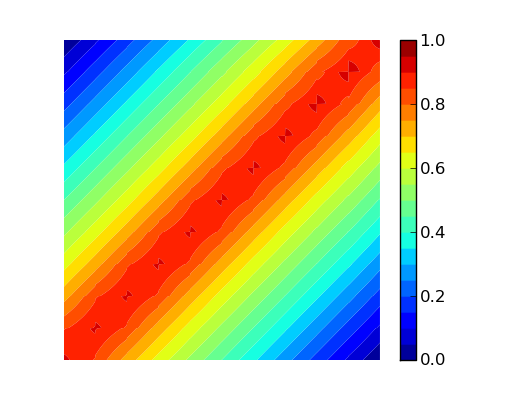}} &
\parbox[c]{\pad\textwidth}{\includegraphics[width=\pad\textwidth]{\pathtofig 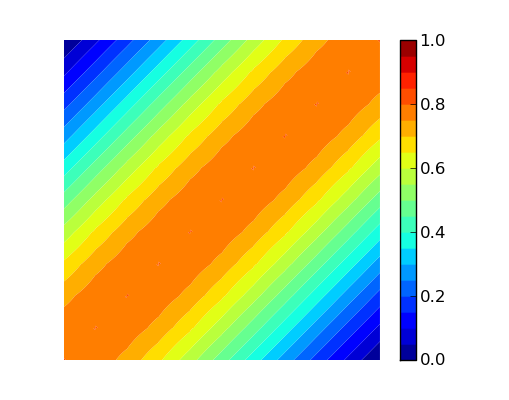}} \\
 & \minmax{0.00}{1.15} & \minmax{0.01}{0.95} & \minmax{0.00}{0.80} \\
\footnotesize{Reference} &
\parbox[c]{\pad\textwidth}{\includegraphics[width=\pad\textwidth]{\pathtofig 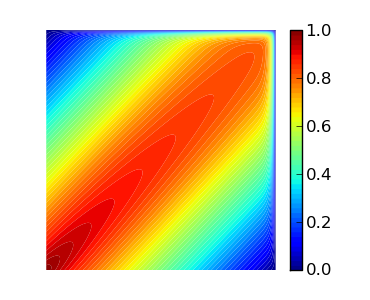}} &
\parbox[c]{\pad\textwidth}{\includegraphics[width=\pad\textwidth]{\pathtofig 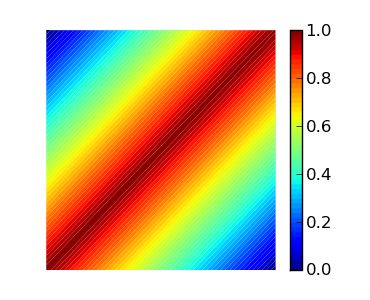}} &
N/A \\
\\
\hline
\end{tabular}
\end{center}
\end{figure}

We observe similar performance of the SN and WN as compared to the
previous problem (\ref{t:legal}). The QON in this case captures the
main features of the solution well. However, while the WN exhibits a
local overshoot in the solution for the $\epsilon=10^{-6}$ case, the
solution value for the QON is underestimated across the
domain. This occurs probably because the test functions are not sufficiently well resolved by our sub-grid discretization, see comments ahead.
\subsubsection*{Advection Skew to the Mesh: Discontinuous Dirichlet Data}
\begin{figure}
\centering
\includegraphics[scale=1.0]{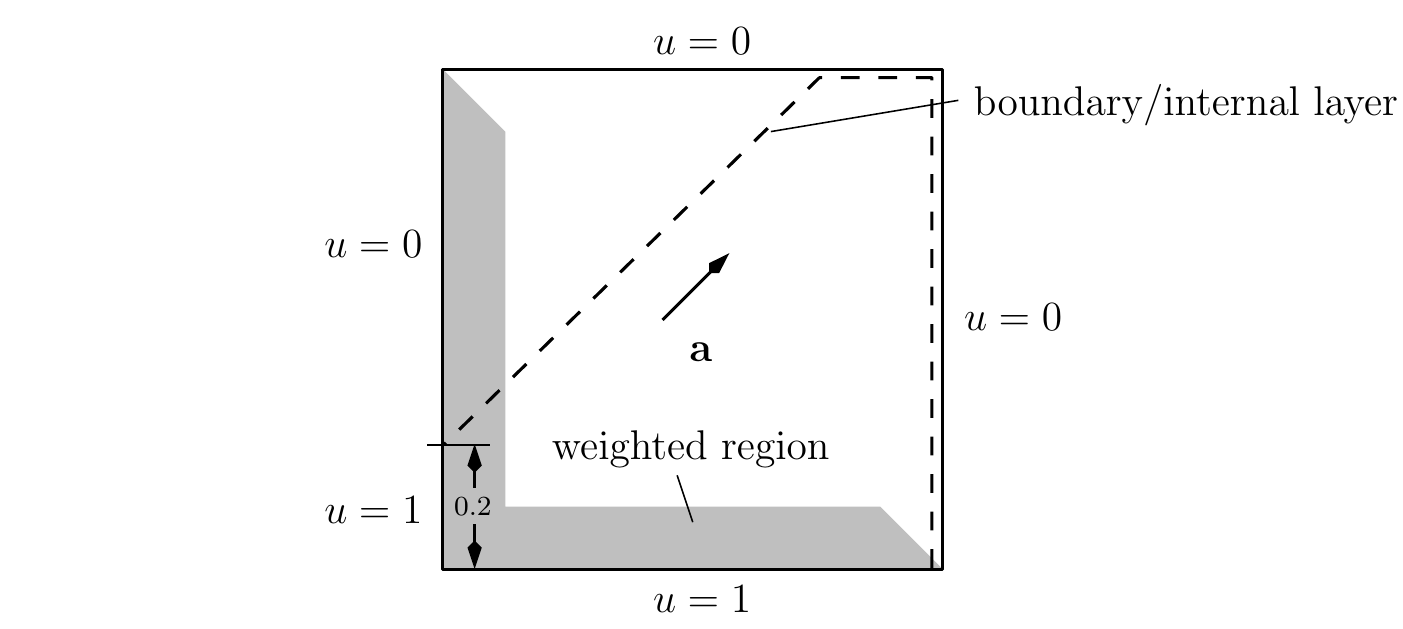}
\caption{Schematics of the advection skew to the mesh problem with discontinuous Dirichlet data. Approximate location of the boundary layer shown as a dashed line and the weighted region is shaded.}
\label{fig:eng}
\end{figure}

\begin{figure}
\caption{Color map plots of $u$ in the advection skew to the mesh problem with discontinuous Dirichlet data. Color range is uniformly set to $[0,1]$ and the minimum and maximum values of the solution are shown underneath each plot.}\label{t:engineering}
\begin{center}
\begin{tabular}{cccc}
\hline
& $\epsilon = 10^{-2}$ 
& $\epsilon = 10^{-4}$ 
& $\epsilon = 10^{-6}$ 
\\
\hline
\\
\footnotesize{SN} &
\parbox[c]{\pad\textwidth}{\includegraphics[width=\pad\textwidth]{\pathtofig 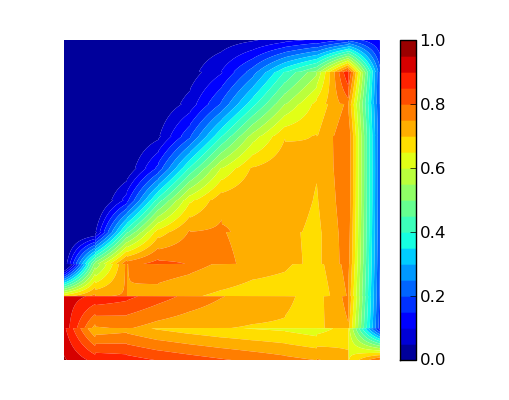}} &
\parbox[c]{\pad\textwidth}{\includegraphics[width=\pad\textwidth]{\pathtofig 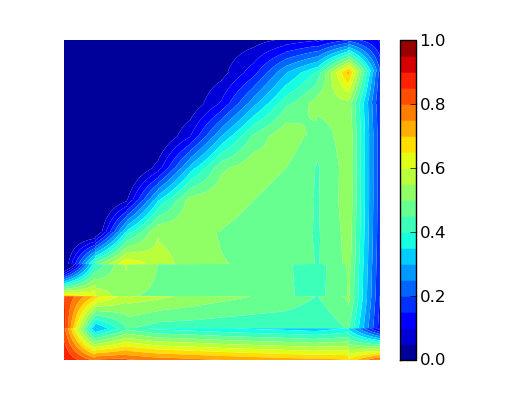}} &
\parbox[c]{\pad\textwidth}{\includegraphics[width=\pad\textwidth]{\pathtofig 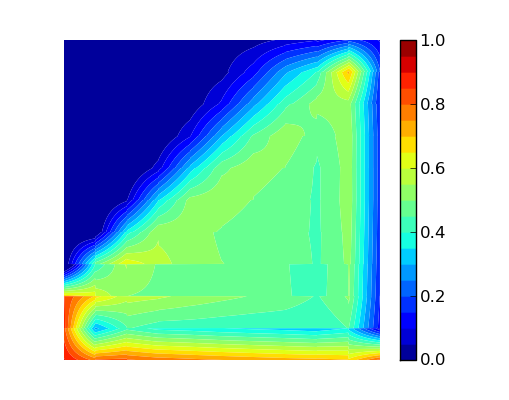}} \\
 & \minmax{-0.13}{0.96} & \minmax{-0.14}{0.91} & \minmax{-0.14}{0.91} \\
\footnotesize{WN} &
\parbox[c]{\pad\textwidth}{\includegraphics[width=\pad\textwidth]{\pathtofig 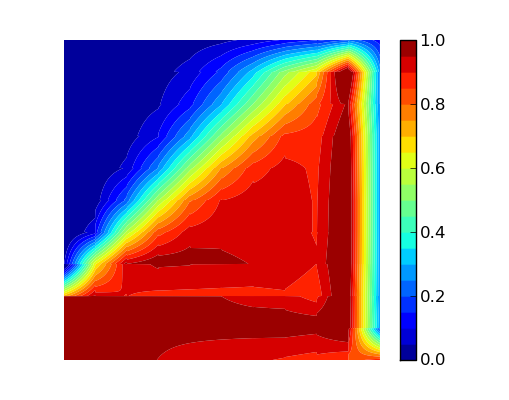}} &
\parbox[c]{\pad\textwidth}{\includegraphics[width=\pad\textwidth]{\pathtofig 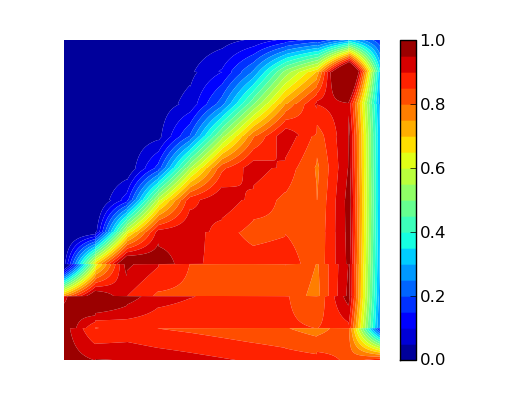}} &
\parbox[c]{\pad\textwidth}{\includegraphics[width=\pad\textwidth]{\pathtofig 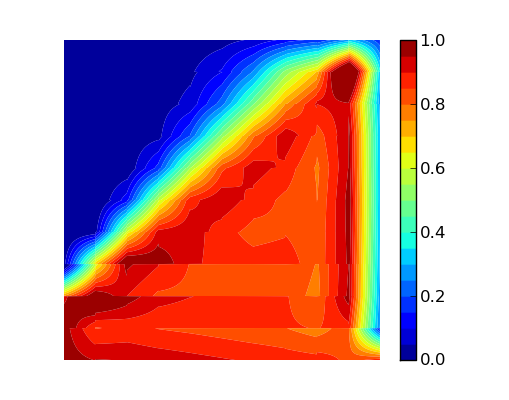}} \\
 & \minmax{-0.09}{1.14} & \minmax{-0.10}{1.26} & \minmax{-0.10}{1.25} \\
\footnotesize{QON} &
\parbox[c]{\pad\textwidth}{\includegraphics[width=\pad\textwidth]{\pathtofig 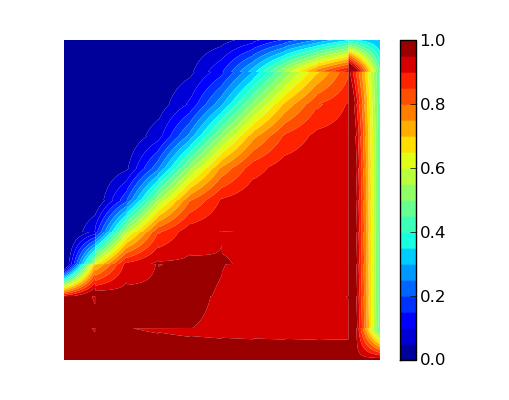}} &
\parbox[c]{\pad\textwidth}{\includegraphics[width=\pad\textwidth]{\pathtofig 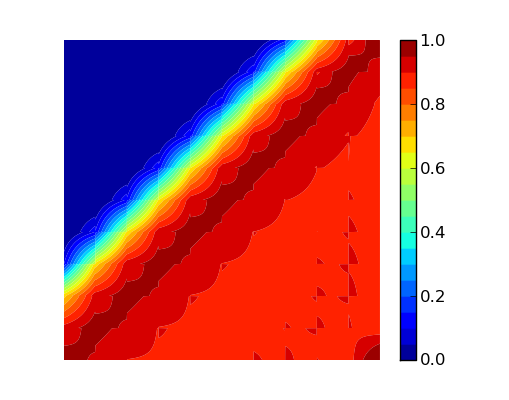}} &
\parbox[c]{\pad\textwidth}{\includegraphics[width=\pad\textwidth]{\pathtofig 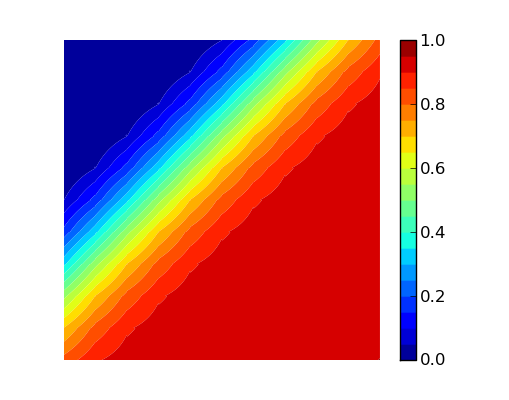}} \\
 & \minmax{-0.12}{1.19} & \minmax{-0.04}{1.04} & \minmax{-0.04}{0.95} \\
\footnotesize{Reference} &
\parbox[c]{\pad\textwidth}{\includegraphics[width=\pad\textwidth]{\pathtofig 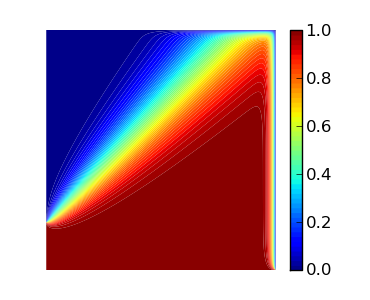}} &
\parbox[c]{\pad\textwidth}{\includegraphics[width=\pad\textwidth]{\pathtofig 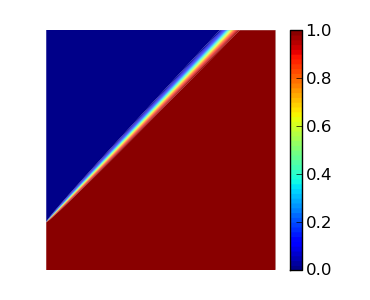}} &
N/A \\
\\
\hline
\end{tabular}
\end{center}
\end{figure}

Our last problem is the standard benchmark test for the
convection-diffusion equation in the engineering literature. The
advection vector is the same as in the previous problem, but the
inflow data is taken to be piecewise constant:
\[
u(x,0) = 1, \quad 
u(y,0) = \begin{cases}
1, &0 \leq y \leq 0.2 \\
0 &0.2 \leq y \leq 1
\end{cases}
\]
The problem schematic is shown in Figure~\ref{fig:eng} and the results
in Figure~\ref{t:engineering}. We emphasize here a feature of QON:
if the mesh does not resolve the boundary layer, its effect is included
but the layer is not resolved. This is in contrast to the SN and WN,
where we always observe a boundary layer the width of one element.
As in the previous problem, the reference solution has been obtained on a fine mesh using
the standard Galerkin FEM.
\subsection{Norm Comparison: Mesh Refinement}
We study the convergence of the different DPG schemes in the advection skew to the mesh problem with continuous Dirichlet data under uniform $h$ and $p$ refinements of the trial space. We compare the error of $u$ and $\msigma$ with respect to the reference solution measured in the $L_2$ norm when $\epsilon = 10^{-2}$. The purpose of this study is to show quantitatively the improved accuracy of the QON solution.

Results from this study are presented in Tables~\ref{tab:trial convergence eps2} and \ref{tab:trial convergence eps4}. Note that in each table, the reference solution comes from the standard Galerkin FEM with a fine mesh of rectangular bilinear elements. The values of $\msigma_{\mathrm{ref}}$ were obtained by a linear interpolation of the coefficients determined from taking finite differences of the solution $u_{\mathrm{ref}}$. In both refinement strategies and for both $u$ and $\msigma$, the QON provides better approximation to the reference solution, particularly when the mesh is coarse.
\begin{table}
\caption{Convergence of $u$ and $\msigma$ under uniform $h$- and $p$-refinement in the advection skew to the mesh problem with continuous Dirichlet data at $\epsilon = 10^{-2}$.}\label{tab:trial convergence eps2}
\begin{center}
\begin{tabular}{cccccccc}
\cline{2-4}\cline{6-8}
 & \multicolumn{3}{c}{$\norm{u_{\mathrm{ref}}-u_n}_{L_2(\Omega)}$} & & \multicolumn{3}{c}{$\norm{\msigma_{\mathrm{ref}}-\msigma_n}_{L_2(\Omega)}$}\\
$N$ & SN & WN & QON & & SN & WN & QON\\
\cline{2-4}\cline{6-8}
5 & 2.55\e{-1} & 1.17\e{-1} & 7.69\e{-2} & & 4.28\e{-2} & 4.15\e{-2} & 3.38\e{-2} \\
10 & 1.76\e{-1} & 6.54\e{-2} & 4.84\e{-2} & & 3.80\e{-2} & 3.68\e{-2} & 2.45\e{-2} \\
25 & 6.26\e{-2} & 2.20\e{-2} & 1.43\e{-2} & & 2.37\e{-2} & 2.33\e{-2} & 9.56\e{-3} \\
50 & 1.64\e{-2} & 7.35\e{-3} & 4.33\e{-3} & & 1.09\e{-2} & 1.08\e{-2} & 3.13\e{-3} \\
100 & 3.24\e{-3} & 2.18\e{-3} & 1.17\e{-3} & & 3.76\e{-3} & 3.74\e{-3} & 8.59\e{-4} \\
\cline{2-4}\cline{6-8}
& & & & & & & \\
& & & & & & & \\
\cline{2-4}\cline{6-8}
 & \multicolumn{3}{c}{$\norm{u_{\mathrm{ref}}-u_n}_{L_2(\Omega)}$} & & \multicolumn{3}{c}{$\norm{\msigma_{\mathrm{ref}}-\msigma_n}_{L_2(\Omega)}$}\\
$p$ & SN & WN & QON & & SN & WN & QON\\
\cline{2-4}\cline{6-8}
1 & 2.55\e{-1} & 1.17\e{-1} & 7.69\e{-2} & & 4.28\e{-2} & 4.15\e{-2} & 3.38\e{-2} \\
2 & 1.49\e{-1} & 5.56\e{-2} & 3.97\e{-2} & & 3.45\e{-2} & 3.29\e{-2} & 2.36\e{-2} \\
3 & 7.00\e{-2} & 2.72\e{-2} & 2.11\e{-2} & & 2.37\e{-2} & 2.28\e{-2} & 1.43\e{-2} \\
4 & 2.58\e{-2} & 1.30\e{-2} & 1.07\e{-2} & & 1.40\e{-2} & 1.37\e{-2} & 7.66\e{-3} \\
\cline{2-4}\cline{6-8}
\end{tabular}
\end{center}
\end{table}
\begin{table}
\caption{Convergence of $u$ and $\msigma$ under uniform $h$- and $p$-refinement in the advection skew to the mesh problem with continuous Dirichlet data at $\epsilon = 10^{-4}$.}
\label{tab:trial convergence eps4}
\begin{center}
\begin{tabular}{cccccccc}
\cline{2-4}\cline{6-8}
 & \multicolumn{3}{c}{$\norm{u_{\mathrm{ref}}-u_n}_{L_2(\Omega)}$} & & \multicolumn{3}{c}{$\norm{\msigma_{\mathrm{ref}}-\msigma_n}_{L_2(\Omega)}$}\\
$N$ & SN & WN & QON & & SN & WN & QON\\
\cline{2-4}\cline{6-8}
5 & 3.73\e{-1} & 1.72\e{-1} & 7.05\e{-2} & & 5.23\e{-3} & 5.23\e{-3} & 5.22\e{-3} \\
10 & 3.60\e{-1} & 1.25\e{-1} & 6.50\e{-2} & & 5.23\e{-3} & 5.22\e{-3} & 5.20\e{-3} \\
25 & 3.51\e{-1} & 9.12\e{-2} & 6.61\e{-2} & & 5.22\e{-3} & 5.21\e{-3} & 5.16\e{-3} \\
50 & 3.45\e{-1} & 7.71\e{-2} & 6.46\e{-2} & & 5.21\e{-3} & 5.18\e{-3} & 5.08\e{-3} \\
100 & 3.36\e{-1} & 6.70\e{-2} & 6.15\e{-2} & & 5.17\e{-3} & 5.13\e{-3} & 4.91\e{-3} \\
\cline{2-4}\cline{6-8}
& & & & & & & \\
& & & & & & & \\
\cline{2-4}\cline{6-8}
 & \multicolumn{3}{c}{$\norm{u_{\mathrm{ref}}-u_n}_{L_2(\Omega)}$} & & \multicolumn{3}{c}{$\norm{\msigma_{\mathrm{ref}}-\msigma_n}_{L_2(\Omega)}$}\\
$p$ & SN & WN & QON & & SN & WN & QON\\
\cline{2-4}\cline{6-8}
1 & 3.73\e{-1} & 1.72\e{-1} & 7.05\e{-2} & & 5.23\e{-3} & 5.23\e{-3} & 5.22\e{-3} \\
2 & 3.56\e{-1} & 1.24\e{-1} & 3.44\e{-2} & & 5.23\e{-3} & 5.22\e{-3} & 5.20\e{-3} \\
3 & 3.52\e{-1} & 1.00\e{-1} & 2.58\e{-2} & & 5.22\e{-3} & 5.20\e{-3} & 5.17\e{-3} \\
4 & 3.49\e{-1} & 9.16\e{-2} & 2.08\e{-2} & & 5.20\e{-3} & 5.18\e{-3} & 5.13\e{-3} \\
\cline{2-4}\cline{6-8}
\end{tabular}
\end{center}
\end{table}
\subsection{The Resolution of the Test Functions Corresponding to the Quasi-optimal Test Space Norm}
The best approximation property of the DPG method hinges on the sufficient resolution of the optimal test functions. In the context of the quasi-optimal norm this is non-trivial because of the singular perturbation character of the auxiliary problem \eqref{eqn:optimal test function problem}. While the subgrid discretization we have proposed is an effort to address this issue, the link between the accuracy of the test functions and the final solution is still unclear. This is partially due to the fact that the auxiliary problem is, to the best of our knowledge, unknown in the literature. Therefore, we make a numerical convergence assessment of the optimal test functions. 

We restrict ourselves to $p$-refinements on the subgrid and measure the approximate relative error of the test functions in the quasi-optimal test space norm $\enorm{\cdot}_{\cW(K)}$ as in
\[
e(\tilde{p}) = \frac{\enorm{\vw^{\tilde{p}+1} - \vw^{\tilde{p}}}_{\cW(K)}}{\enorm{\vw^{\tilde{p}+1}}_{\cW(K)}}
\]
This we carry out on a single square element $K$ of the ten by ten mesh used earlier in the advection skew to the mesh problems. Because the norm $\enorm{\cdot}_{\cW(K)}$ is the usual energy norm for the auxiliary problems, its computation reduces to linear algebra. Recall that the columns of the matrix $\mU_{\mathrm{opt}}$ represent the optimal test functions. Therefore the energy norm squared of the optimal test functions $\vw_i$ is given by 
\[
\enorm{\vw_i}^2_{\cW(K)} = \mU_{\mathrm{opt}}(:,i)^T \mK_{\mathrm{opt}} \mU_{\mathrm{opt}}(:,i) 
= \mU_{\mathrm{opt}}(:,i)^T\mF_{\mathrm{opt}}(:,i)
= \mK_{\ell}(i,i)
\]
where we have adopted the \verb|MATLAB| notation for indexing a whole row or column of a matrix.
In other words, the energy norm of the optimal test functions may be read directly from the diagonal of $\mK_{\ell}$.

Notice that when $p=1$, there are 28 trial functions and optimal test functions. To simplify the representation of the results, we separate the optimal test functions into four groups associated to the trial functions for $\msigma$, $u$, $\hat{u}$, and $\hat{\sigma}_n$ and compute the average of the errors for each group. We present the results in Tables~\ref{tab:test function convergence eps2} and \ref{tab:test function convergence eps4} for $\varepsilon = 10^{-2}$ and $\varepsilon = 10^{-4}$ where we have used $\bar{e}$ to indicate an average error and subscripts to refer to each group. We have also included the $L_2$ error with respect to the reference solution of $u$ in the table. The results show that while the optimal test functions corresponding to the field variables $\msigma$ and $u$ are seemingly easy to approximate, those corresponding to the interface variables $\hat{\sigma}_n$ and $\hat{u}$ are problematic. As a matter of fact, the relative error level of the optimal test functions corresponding to the trace $\hat{u}$ barely decreases below the discretization error level for the field varible $u$ in the trial space as $\Delta p$ increases. This slow convergence explains the lull in the $h$-convergence in Table \ref{tab:trial convergence eps4} as well as the degradation of the global accuracy when $\epsilon = 10^{-6}$ in Figures \ref{t:legal} and \ref{t:engineering}. These phenomena are manifestations of the current practical limitations in the utilization of the quasi-optimal norm. 
\begin{table}
\caption{Convergence of the optimal test functions in the advection skew to the mesh problems when $\epsilon= 10^{-2}$, $N=10$ and $p=1$.}
\label{tab:test function convergence eps2}
\begin{center}
\begin{tabular}{ccccccc}
\cline{2-5}\cline{7-7}
 & \multicolumn{4}{c}{average relative approximate errors} & & \\
$\Delta p$ & $\bar{e}_{\msigma}$ & $\bar{e}_{u}$ & $\bar{e}_{\hat{u}}$ & $\bar{e}_{\hat{\sigma}_{n}}$ & & $\norm{u_{\mathrm{ref}}-u_n}_{L_2(\Omega)}$ \\
\cline{2-5}\cline{7-7}
1 & 2.16\e{-3} & 4.85\e{-3} & 9.42\e{-2} & 1.13\e{-3} & & 4.62\e{-2} \\
2 & 3.28\e{-4} & 7.38\e{-4} & 1.42\e{-2} & 1.67\e{-4} & & 4.79\e{-2} \\
3 & 3.08\e{-5} & 6.74\e{-5} & 1.77\e{-3} & 1.57\e{-5} & & 4.84\e{-2} \\
4 & 2.52\e{-6} & 5.34\e{-6} & 2.15\e{-4} & 1.30\e{-6} & & 4.85\e{-2} \\
5 & 4.22\e{-7} & 8.70\e{-7} & 4.23\e{-5} & 2.18\e{-7} & & 4.85\e{-2} \\
6 & 1.08\e{-7} & 2.20\e{-7} & 1.12\e{-5} & 5.58\e{-8} & & 4.85\e{-2} \\
7 & 3.44\e{-8} & 6.92\e{-8} & 3.59\e{-6} & 1.77\e{-8} & & 4.85\e{-2} \\
8 & 1.27\e{-8} & 2.53\e{-8} & 1.32\e{-6} & 6.54\e{-9} & & 4.85\e{-2} \\
9 & 5.22\e{-9} & 1.04\e{-8} & 5.45\e{-7} & 2.69\e{-9} & & 4.85\e{-2} \\
\cline{2-5}\cline{7-7}
\end{tabular}
\end{center}
\end{table}
\begin{table}
\caption{Convergence of the optimal test functions in the advection skew to the mesh problems when $\epsilon= 10^{-4}$, $N=10$ and $p=1$.}
\label{tab:test function convergence eps4}
\begin{center}
\begin{tabular}{ccccccc}
\cline{2-5}\cline{7-7}
 & \multicolumn{4}{c}{average relative approximate errors} & & \\
$\Delta p$ & $\bar{e}_{\msigma}$ & $\bar{e}_{u}$ & $\bar{e}_{\hat{u}}$ & $\bar{e}_{\hat{\sigma}_{n}}$ & & $\norm{u_{\mathrm{ref}}-u_n}_{L_2(\Omega)}$ \\
\cline{2-5}\cline{7-7}
1 & 2.07\e{-4} & 5.40\e{-3} & 4.54\e{-1} & 5.36\e{-3} & & 9.15\e{-2} \\
2 & 1.51\e{-5} & 3.00\e{-5} & 3.25\e{-1} & 7.65\e{-4} & & 7.22\e{-2} \\
3 & 8.17\e{-6} & 1.68\e{-5} & 1.48\e{-1} & 3.99\e{-4} & & 6.50\e{-2} \\
4 & 3.10\e{-6} & 6.22\e{-6} & 4.67\e{-2} & 1.50\e{-4} & & 6.32\e{-2} \\
5 & 1.17\e{-6} & 2.18\e{-6} & 1.43\e{-2} & 5.75\e{-5} & & 6.29\e{-2} \\
6 & 5.36\e{-7} & 9.49\e{-7} & 5.77\e{-3} & 2.67\e{-5} & & 6.29\e{-2} \\
7 & 2.97\e{-7} & 5.02\e{-7} & 3.56\e{-3} & 1.49\e{-5} & & 6.29\e{-2} \\
8 & 1.85\e{-7} & 3.09\e{-7} & 2.76\e{-3} & 9.24\e{-6} & & 6.30\e{-2} \\
9 & 1.24\e{-7} & 2.04\e{-7} & 2.39\e{-3} & 6.15\e{-6} & & 6.31\e{-2} \\
\cline{2-5}\cline{7-7}
\end{tabular}
\end{center}
\end{table}

\newpage
\section{Concluding Remarks} \label{sec:conclusions}
We have studied the robustness of the discontinuous Petrov-Galerkin method with optimal test functions when applied to convection-diffusion problems. In particular, we have introduced a DPG scheme based on the quasi-optimal test norm and compared it with schemes based on other test space norms. Our numerical tests reveal that the quasi-optimal norm leads to improved coarse mesh accuracy as compared with the conventional DPG test space norms provided that the associated test functions are sufficiently well resolved. We have addressed the resolution of the optimal test functions by studying their convergence numerically. We have been able to find a compromise between accuracy and stability of the final DPG solution by combining a carefully designed element sub grid with numerical regularization of the test functions.

The DPG method with optimal test functions originates from the functional analytic foundation of the finite element method, and is perhaps still far from being adopted into engineering practice. We have made an attempt to popularize the methodology by contributing a holistic approach which begins in the abstract framework and ends in detailed explanation of the algorithm from the implementation viewpoint.

\bibliographystyle{acm}
\bibliography{Refer}

\end{document}